\documentclass[10pt,leqno]{amsart}
\usepackage{comment}
\usepackage{graphicx}
\usepackage{amssymb,amsthm,amsmath}

\newtheorem{theorem}{Theorem}[]
\newtheorem{definition}[theorem]{Definition}
\newtheorem{example}[theorem]{Example}
\newtheorem{lemma}[theorem]{Lemma}
\newtheorem{proposition}[theorem]{Proposition}

\newtheorem{remark}[theorem]{Remark}
\usepackage{lipsum}

\newcommand{\R}[1]{{\mathbb R}^{#1}}

\newcommand{\1}{\mathbf{1}}

\newenvironment{pf}{\noindent{\it Proof.} }{\;\begin{flushright}\fbox{}\end{flushright}}

\usepackage{amsmath}
\usepackage{amsfonts}
\usepackage{amssymb}
\usepackage{bbm}
\usepackage{graphicx}
\usepackage{multicol}
\usepackage{fancybox}
\usepackage[T1]{fontenc}
\usepackage{frcursive}
\usepackage{calligra}
\usepackage{marginnote}
\usepackage{dsfont}
\usepackage{float}
\usepackage{physics}
\usepackage{subfigure}
\usepackage{dirtytalk}
\usepackage{csquotes}
\usepackage{hyperref}
\usepackage{mathrsfs}
\usepackage{arydshln}
\usepackage{lineno}
\makeatletter
\renewcommand*\env@matrix[1][*\c@MaxMatrixCols c]{%
  \hskip -\arraycolsep
  \let\@ifnextchar\new@ifnextchar
  \array{#1}}
\makeatother

\DeclareMathAlphabet{\mathpzc}{OT1}{pzc}{m}{it}

\newcommand{\prob}     {\mathbb{P}}        
\newcommand{\E}     {\mathbb{E}}        
\newcommand{\XX}{{\bf X}}
\newcommand{\YY}{{\bf Y}}
\newcommand{\xx}{{\bf x}}
\newcommand{\yy}{{\bf y}}

\newcommand{\WW}{{\bf W}}
\newcommand{\ZZ}{{\bf Z}}
\newcommand{\zz}{{\bf z}}

\newcommand{\A}{{\bf A}}

\newcommand{\TT}    {\mathcal{T}}

\newcommand{\FF}    {\mathcal{F}}

\newcommand{\dist} {\mathpzc{ dist}}
\newcommand{\interior}{\mathpzc{int}}
\newcommand{\boundary}{\partial}

\newcommand{\diag}{\text{\textbf{diag}}}
\newcommand{\Lip}{\textbf{Lip}}
\newcommand{\Law}{\textbf{Law}}
\newcommand{\Was}{\textbf{Was}}
\newcommand{\Couplings} {\textbf{Couplings}}
\newcommand\eqLaw{\mathrel{\stackrel{\makebox[0pt]{\mbox{\normalfont\tiny Law}}}{=}}}

\newcommand{\inner}[2] { \left \langle #1 , #2  \right \rangle }
\newcommand{\bG} {{\bf G}}

\newtheorem{assumptions}{Assumptions}

\newcommand{\simp}[1]{\mathbb{S}^{#1}}

\newcommand{\proj}{\Pi_\mathcal{T}}

\begin{document}

\title[Interacting replicator dynamics]{Persistence and neutrality in interacting replicator dynamics \thanks{This work was funded by project FONDECYT 1200925 \emph{The emergence of ecologies through metabolic cooperation and recursive organization}, Centro de Modelamiento Matem\'atico (CMM),
Grant FB210005, BASAL funds for centers of excellence from ANID-Chile, Exploration-ANID 13220168, \emph{Biological and quantum Open System Dynamics: evolution, innovation and mathematical foundations}, FONDECYT Iniciaci\'on project number 11240158-2024 \emph{Adaptive behavior in stochastic population dynamics and non-linear Markov processes in ecoevolutionary modeling}, and FONDECYT Iniciaci\'on 11200436, \emph{Excitation and inhibition balance as a dynamical process.}, and \textit{``Programa de Inserción Acad\'emica 2024 Vicerrector\'ia Acad\'emica y Prorrectoría de la Pontificia Universidad Cat\'olica de Chile''}}
}




\author[Videla]{Leonardo Videla$^1$}
\address{$^1$Departamento de Matem\'atica y Ciencia de la Computaci\'on, Facultad de Ciencia, Universidad de Santiago de Chile, Las Sophoras 173, Estaci\'on Central, Santiago, Chile}\email{leonardo.videla@usach.cl}

\author[Tejo]{Mauricio Tejo$^2$}
\address{$^2$Instituto de Estad\'istica, Universidad de Valpara\'iso, Valpara\'iso, Chile}
\email{mauricio.tejo@uv.cl}

\author[Qui\~ninao]{Crist\'obal Qui\~ninao$^3$}
\address{$^3$Facultad de Ciencias Biol\'ogicas, Pontificia Universidad Cat\'olica de Chile, Santiago, Chile}
\email{cquininao@uc.cl}
\email{pmarquet@puc.bio.cl}

\author[Marquet]{Pablo A. Marquet$^{3,4}$}
\address{$^4$Centro de Modelamiento Matem\'atico (CMM), Universidad de Chile-IRL 2807 CNRS Beauchef 851, Santiago, Chile}
\address{$^4$The Santa Fe Institute, 1399 Hyde Park Road, Santa Fe NM 87501, USA}
\address{$^4$Instituto de Sistemas Complejos de Valpara\'iso, Subida Artiller\'ia 470, Valpara\'iso, Chile}

\author[Rebolledo]{Rolando Rebolledo$^5$}
\address{$^5$Instituto de Ingenier\'ia Matem\'atica, Universidad de Valpara\'iso}
\email{rolando.rebolledo@uv.cl}

\date{Received: date / Accepted: date}

\begin{abstract}
We study the large-time behavior of an ensemble of entities obeying replicator-like stochastic dynamics with mean-field interactions as a model for a primordial ecology. We prove the propagation-of-chaos property and establish conditions for the strong persistence of the $N$-replicator system and the existence of invariant distributions for a class of associated McKean-Vlasov dynamics. In particular, our results show that, unlike typical models of neutral ecology, fitness equivalence does not need to be assumed but emerges as a condition for the persistence of the system. Further, neutrality is associated with a unique Dirichlet invariant probability measure. We illustrate our findings with some simple case studies, provide numerical results, and discuss our conclusions in the light of Neutral Theory in ecology.


\end{abstract}

\maketitle

\section{Introduction} \label{sec:introduction}

The emergence of life corresponds to the emergence of an ecological system of interacting self-replicating entities.  Several characteristics of early life have endured through time and are shared by all forms of life today, from the near universality of the genetic code to the intermediate metabolism \cite{koonin2009origin,smith2016origin}. It is unlikely that these biological universalities represent only "frozen accidents" linked to a universal common ancestry, but may emerge as a consequence of some fundamental principles associated with information and thermodynamics that affect the robustness and evolvability of biological systems \cite{smith2016origin,walker2017matter}. In a similar vein, it is expected that early ecological systems faced important challenges to their persistence in fluctuating environments, where the number of co-existing entities increases forming diverse communities, and that some fundamental principles may also be invoked to understand their persistence.  In this work, we are primarily interested in the \emph{principle of neutrality}, which is mathematically expressed in a particular invariant probability distribution associated with the process of abundances of different types of entities (e.g. species) in an ecological system (e.g. community) where birth and death rates are linear, and thus entities can be considered as weakly interacting (e.g., a neutral ecology where species are equivalent in fitness, \cite{hubbell2001}) due to the operation of fitness equalizing mechanisms associated to tradeoffs and incompatible optima. These mechanisms render neutrality as a good approximation, even if entities weakly interact, as envisioned in the Red Queen and in the Neutral Theory of Ecology \cite{van1973new,hubbell2001,hubbell2005neutral,hubbell2006neutral}. In particular, using the neutrality approximation in a stochastic framework, we showed in~\cite{marquet2017proportional} that a one-dimensional diffusion approximation for the frequency or proportional abundance of species in a community admits a Beta invariant distribution (the Proportional Species Abundance Distribution), the same as found, for instance, for gene frequencies in partially isolated populations \cite{wright1931evolution,kolmogorov1935}, which generalizes to a Dirichlet in the multidimensional case.   

To model an early ecological system where to shed light on the importance of neutrality in early ecologies,  we use a quintessential model for replicating entities, which has been applied to the modeling of polynucleotides in a container, the dynamics of gene frequencies, selection (Price equation), relative densities of interacting populations (Lotka-Volterra model), the frequency of different strategies in a population (game theoretical models), and adaptive dynamics \cite{taylorjonker1978,schuster1983replicator,hofbauerbook,Page:2002aa}. One of the most basic and general equations for all these dynamics is the deterministic equation introduced by Taylor and Jonker \cite{taylorjonker1978} and dubbed \emph{replicator equation} by Schuster and Sigmund \cite{taylorjonker1978,schuster1983replicator}. 

Thus, in this article, any given entity in our primordial ecological system will be represented by a replicator, and we will follow its proportional abundance through time. Now, since it can be argued that in the primordial soup, many entities shared a common habitat in this article we will extend the simple replicator dynamics to the setup of a large number of self-reproducing entities (community), constituting in this way an interactive set. In the setting we have chosen, the intraspecific dynamics of each singular replicator will be affected, in the form of a mean-field effect, by the dynamics of the other replicators in the system. 

\emph{Is neutrality still plausible to emerge in this larger setting? and if so, What is its role?} The sections that follow are devoted to shed light on these questions.

\subsection{Preliminaries}

The main motivation of Taylor and Jonker's replicator model (\cite{taylorjonker1978}) was to unify many game-theoretic models that pursue the study of evolutionary stable strategies (ESS) in the idealized dynamics of animal conflicts as introduced by Maynard Smith and Price \cite{smith1973logic}. Recall that Taylor and Jonker's replicator is a dynamical system $(\xx(t): t \ge 0)$ taking values on the $(d-1)$-dimensional simplex:
\begin{align*}
    \simp{d}:= \{\xx \in \mathbb{R}^{d}: x_i \ge 0, \sum_{i=1}^{d} x_i=1\},
\end{align*}
obeying the ordinary differential equation:
\begin{align} \label{eq:original_replicator}
\dfrac{d\xx}{dt}= \xx \circ (A\xx- \inner{\xx}{A\xx} \mathbf{1}),
\end{align}
for some initial condition $\xx_0 \in \simp{d}$. Here, $\circ$ stands for the Hadamard product (entry-wise), $A$ is a suitable $d \times d$ real matrix termed the payoff matrix, and $\mathbf{1}$ stands for the vector of dimensions $d$ with all its entries equal to $1$. Equation \eqref{eq:original_replicator} is intended to model the dynamics of the proportions $(x_i: i=1, \ldots, d)$ of individuals playing contesting strategies in (the continuous-time approximation of) a multi-round, two-player game: if each player can choose among strategies $1, 2, \ldots, d$, the entry $A_{ij}$ is interpreted as the payoff from using strategy $i$ against an individual using strategy $j$. The term $- \inner{\xx}{A\xx} \mathbf{1}$ modulates the fitness variation in such a way that the fitness of a player of type $i$ increases if the payoff is greater than the average payoff of the community (see \cite{hofbauerbook,Hofbauer:1981ab,Hofbauer:1981aa} for a thorough study of these systems in the context of evolutionary game theory). 

Many generalizations of the original replicator dynamics can be considered. For example, 
we can obtain the generalized replicator dynamics:
\begin{align*}
    \dfrac{d\xx}{dt}= \xx \circ (F(\xx)- \inner{\xx} {F(\xx)} \mathbf{1}),
\end{align*}
for some suitable instantaneous fitness function $F: \simp{d} \mapsto \mathbb{R}^{d}$. In this form, replicator dynamics has become a general model for the evolution of populations under frequency-dependent selection, ``flexible enough to cover a great deal of evolution models, suggesting a unifying view of replicator selection from the primordial soup up to animal societies'' (\cite{schuster1983replicator}). The mathematical analysis of general replicator equations has been extensive (\cite{hofbauerbook1988}, \cite{hofbauerbook}, \cite{Hofbauer2003b}, \cite{weibull1995evolutionary}, \cite{hauert2002volunteering}, \cite{cressman2003evolutionary}, \cite{nowak2004evolutionary}, \cite{ohtsuki2006replicator}, \cite{cressman2014replicator}), \cite{cooney2019replicator},  testifying to its fundamental role in the development of evolutionary game theory. 


Stochastic versions of the replicator equation have a long tradition. From the initial efforts of Foster and Young (\cite{fosteryoung90}) and Fudenberg and Harris (\cite{fudenbergharris92}) in exploring the replicator equation in a stochastic differential equation form (see also \cite{amir1998scale}, \cite{saito1997note}), the analysis of stochastic replicator equations have become increasingly complex, going beyond the case of only two pure strategies (\cite{imhof2005}) to replicator dynamics with Stratonovich-type perturbations (\cite{khasminskii2006replicator}), to the analysis of stochastic persistence (\cite{benaim2008robust}), variations in game dynamics with a potential function (\cite{avrachenkov2019metastability}) to the inclusion of revision protocols such as imitation (\cite{sandholm2010population}, \cite{mertikopoulos2016imitation}). 

It must be said that the interaction between the main system and the noise can be interpreted in many different ways. For example, Fudenberg and Harris's version of the stochastic replicator can be obtained from the following considerations regarding the instantaneous fitness (see \cite{fudenbergharris92}). Let $\mathbb{R}^{d}_{++}$ be the cone of strictly positive real $d$-tuples, let $\YY_t= (Y^{(1)}_t, Y^{(2)}_t,\ldots ,Y^{(d)}_t )^\top \in \mathbb{R}^{d}_{++}$ be the vector of abundances at time $t \ge 0$ of the $d$ types in a focal population, and let:
\begin{align*}
\XX_t=\dfrac{\YY_t}{\sum_{i=1}^d Y^{(i)}_t} \in \simp{d}, 
\end{align*}
be the relative proportions at time $t$ of these types. We prescribe a stochastic evolution of the abundances through the equation. 
\begin{align}\label{eq:abundance_dynamics}
    d\YY_t= \YY_t \circ (A \XX_t dt + \Sigma d \WW_t),
\end{align}
where $\Sigma : = \diag(\sigma_1, \sigma_2, \ldots, \sigma_d)$ is a diagonal non-singular matrix and $\WW= (W_1, \ldots, W_d)$ is a standard $d$-dimensional Wiener process defined on a probability space $(\Omega, \FF, \prob)$ and adapted to a right-continuous and augmented filtration $\mathbb{F}:= (\FF_t: t \geq 0)$. The rationale behind this evolution is easy to grasp: the instantaneous fitness of type $i$ is the payoff associated with the strategy of playing $i$ at the density level $\XX$, \emph{plus} a random perturbation. In \eqref{eq:abundance_dynamics} the stochastic term is generally interpreted as environmental fluctuations due, for example, to high-frequency variations of the abiotic or biotic factors that underlie the ecological dynamics. Depending on the time scale of the phenomenon we are interested in, these variations can come from different sources: temperature fluctuations, changes in resource availability, among others. 


From equation \eqref{eq:abundance_dynamics}, it is easy to show that the dynamics of the proportion obeys the SDE: 
\begin{align} \label{eq:stochastic_replicator_1}
    d\XX_t =  \XX_t \circ (\tilde A \XX_t - \inner{\XX_t}{\tilde A \XX_t}\mathbf{1} ) dt + \sum_{i=1}^d \XX_t \circ \left ( \Sigma e_i - \inner{\XX_t}{\Sigma e_i} \mathbf{1}\right) d W_t^{(i)}, 
\end{align}
where $\tilde A= A- \Sigma \Sigma^\top$ and $e_i$ is the $i$-th vector of the canonical basis of $\mathbb{R}^d$. 


Hofbauer and Imhof \cite{hofbauerimhof2009} studied  equation \eqref{eq:stochastic_replicator_1}, with a focus on the time-averages of the dynamics. Among other results, the authors show that under some conditions on the payoff matrix, the process is positive recurrent and its transition probability function admits an invariant Dirichlet distribution. Interestingly, this holds when the sum of the off-diagonal payoffs and the diagonal obey a simple linear relation (this is their Corollary 3.12). 
This symmetry can be regarded as equivalent to a type of neutrality, as described in the first paragraphs of this introduction,  and in fact, induces the emergence of a neutral ecology in replicator dynamics, which provides a remarkable connection to the existing neutral theory in community  ecology \cite{hubbell2001,marquet2017proportional}.  
Thus, in this contribution we will use the replicator equation in a stochastic framework as a vehicle to understand the emergence of primordial ecologies and test for the emergence and generality of the fitness equivalence principle outlined in \cite{hofbauerimhof2009,marquet2017proportional}, in the context of a community of interacting replicators, as explained next. 

\subsection{Our model: interacting replicators}
For measurable functions $F: \simp{d} \mapsto \mathbb{R}^d$, write:
\begin{align*}
    \proj F (\xx) := \xx \circ (F (\xx) - \inner{\xx}{F(\xx)}\mathbf{1}),
\end{align*}
and extend the definition to matrix-valued functions by tensorisation, namely: if $\mathbf{F}= (F_1, \ldots, F_m)$ where $F_i: \simp{d} \mapsto \mathbb{R}^d$ is bounded and measurable, then $\proj \mathbf{F}$ is the $d \times m$ matrix whose $j$-th column is $\proj F_j (\xx)$. 
With this, equation \eqref{eq:stochastic_replicator_1} can be written as:
\begin{align}\label{eq:eq00}
    d\XX_t = \proj \Phi (\XX_t) dt + \proj \Psi (\XX_t) d \WW_t,
\end{align}
for a suitable vector field $\Phi$ and a suitable matrix-valued function $\Psi$; here $\WW$ is a standard $d$-dimensional Brownian motion on a filtered probability space. This is the basic replicator equation we will be working with in this contribution. 

Now, to capture the interaction of one replicator with other  of the same type, we assume that a smooth function $\Upsilon: \simp{d} \times \simp{d} \mapsto \mathbb{R}^d$ is given. The notation $\proj \Upsilon (\xx, \yy)$ stands for the action of $\proj$ on the first variable. Let be given a natural number $N$, and consider the $(\simp{d})^N$-valued stochastic process defined via 
\begin{align} \label{eq:mean_field_SDE}
  \nonumber  d\XX^{(N; i)}_t &= \proj \Phi (\XX^{(N;i)}_t) dt 
    \\
    &+ \dfrac{1}{N} \sum_{j=1}^N \proj \Upsilon (\XX^{(N;i)}_t, \XX^{(N;j)}_t) dt + \proj \Psi (\XX^{(N;i)}_t) d\WW^{(i)}_t,
\end{align}
 for $i=1, 2, \ldots, N$. Here, $(\WW^{(i)}, i=1, 2, \ldots, N)$ is a family of independent $\mathbb{F}$-adapted, $d$-dimensional standard Wiener processes carrying the local noises (i.e., noise associated to each species represented by a single replicator). 

 To fix ideas, consider the case $d=2$, $\Phi$ linear and $\Psi$ a constant matrix, and assume, for instance, that $\Upsilon$ has the form:

\begin{align*}
\Upsilon (\xx, \yy) = \delta y_1 \mathbf{E} \xx,
\qquad
\mathbf{E}= \Bigg ( \begin{matrix}
    0 & 1 \\
    -1 & 0 
\end{matrix} \Bigg)
\end{align*}
where $\delta$ is a small parameter. We will return to this example in Section \ref{section:invariant}, and so it seems just right to use it to illustrate the meaning of the equations above. Indeed, in this case, we can provide an intuitive interpretation of equation~\eqref{eq:mean_field_SDE}: the payoff obtained from using the first strategy against a player using the second one increases proportionally to the overall average frequency of the usage of the first strategy amongst the family of replicators. Thus, at least intuitively, in this system of interacting replicators, the first strategy becomes more resilient, and thus one should expect higher average usage of the first strategy as compared to the situation with no interaction. Of course, the overall effect must take account of the particular form of the intrinsic payoff matrix (the diagonal elements could restrain a net bias to the use of the first strategy) and the noise part (if it is large, randomness may dominate the dynamics), but as we will see later, under appropriate assumptions, the intuition above is not far from the actual long-term picture.  

\begin{remark}
In this article, we are considering the case of interacting replicators of just one type, which could be akin to traditional species that are either phylogenetically and functionally close and thus share a similar niche, trophic level, and interactions within communities, as envisioned in the principle of niche conservatism \cite{wiens2010niche}. As we will see below, this will translate in that, for large $N$, the individual replicators behave as statistically identical entities. The multi-type model (which would amount to mean-field interactions of $N_1$ replicators of type $1$, $N_2$ replicators of type $2$, etc.; see for instance models of this type associated to neurosciences~\cite{faugeras2009constructive}, social sciences~\cite{contucci2008phase} or statistical mechanics~\cite{collet2014macroscopic}) is technically more involved and represents a reasonable second step in our research program. 
\end{remark}

Systems analogous to equation \eqref{eq:mean_field_SDE} are generically referred to as \emph{mean-field interacting particle systems}, and one of the main problems associated with them is the determination of conditions under which, in the thermodynamic limit $N \to \infty$, the assembly of particles (in our case, replicators) behave as if they were independent. This is the \emph{propagation of chaos property}, whose precise content will be given later. For the time being, assume that we are interested in the study of the system \eqref{eq:mean_field_SDE} on a finite time horizon and for \emph{large values of} $N$. Write $\mu_t^{(N)}$ for the empirical probability measure induced on the simplex by $\vec \XX^{(N)}_t:= (\XX^{(N;i)}_t, \XX^{(N;2)}_t, \ldots, \XX^{(N;N)}_t)^\top$, namely:
\begin{align*}
    \mu_t^{(N)} (A) := \dfrac{1}{N} \sum_{i=1}^N \delta_{\XX^{(N;i)}_t} (A).
\end{align*}
With a similar notation, $\mu_{[0, t]}^{(N)} (A)$ stands for the random probability measure on $\mathcal{C}([0, t]; \simp{d})$ given by:
\begin{align}
\label{eq:random_pm}
\mu^{(N)} (A) =\mu_{[0, t]}^{(N)} (A) := \dfrac{1}{N} \sum_{i=1}^N \delta_{\XX^{(N;i)}_{[0, t]}} (A).  
\end{align}
The following convention will be used in the sequel. For an arbitrary metric space $E$ and a measurable function  $G: \simp{d} \times \simp{d} \mapsto E$, we will write:
\begin{align*}
    G(\xx, \mu):= \int_{\simp{d}} G (\xx, \yy)\mu (d\yy),
\end{align*}
and thus, we can regard $G$ as a map $\simp{d} \times \mathcal{P}(\simp{d}) \mapsto E$. With this convention, equation \eqref{eq:mean_field_SDE} can be written as: 
\begin{align}\label{eq:MFSDE} 
    d\XX^{(N;i)}_t& = \proj \Phi (\XX^{(N;i)}_t) dt +  \proj \Upsilon (\XX^{(N;i)}_t, \mu^{(N)}_{t}) dt  \\
    & \qquad \qquad + \proj \Psi (X_t^{(N; i)}) d\WW^{(i)}_t,\qquad i=1, 2, \ldots, N. \nonumber
\end{align}
Assume that we know that for each $t$ in a finite time horizon, the random probability measures $(\mu_t^{(N)}: N \ge 1)$ converge (almost surely in the weak topology of probability measures) to a certain deterministic probability measure, $\mu_t$ say. Then, at least intuitively, we should have that for large values of $N$, the dynamics of a typical particle should be governed by:
\begin{align} \label{eq:mckean_vlasov_0}
    d\Xi_t = \proj \Phi (\Xi_t) dt + \proj \Upsilon (\Xi_t, \mu_t)dt + \proj \Psi (\Xi_t) d\WW_t, 
\end{align}
and, following a Glivenko-Cantelli-like reasoning, it should be apparent that $\mu_t$ is nothing but the law of $\Xi_t$. Equation \eqref{eq:mckean_vlasov_0} is of the McKean-Vlasov type, and its solutions (whenever they exist) induce a paradigmatic example of a \emph{non-linear Markov process}. Here, the term non-linear has nothing to do with the algebraic degree of the unknown $\Xi$ in \eqref{eq:mckean_vlasov_0} but rather with the fact that the random dynamics depends on the very distribution of the Markov process; in turn, it translates into a non-linearity at the level of the evolution equation associated to the dynamics. There is a rich literature on McKean-Vlasov SDE, with the connection with nonlinear PDE being an important research topic (see the classical \cite{mckean66} for the motivations underlying this kind of process).

The aim of this article is to study the link between the mean-field system \eqref{eq:mean_field_SDE} and the McKean-Vlasov equation \eqref{eq:mckean_vlasov_0}; to study the persistence (in a sense to be defined below) for the mean-field system and to prove existence and uniqueness of non-trivial invariant measures for the McKean-Vlasov equation in some simple cases. All these issues are treated with an eye on the problem posed at the opening of the introduction, namely the emergence and role of neutrality in large interacting sets of biological entities. The organization of the paper reflects closely this sequence of problems. In Section \ref{sec:framework} we present all the mathematical preliminaries that allow us to put the models on firm grounds, and we describe our assumptions and main results. Then, we prove the analytical properties of the processes and provide a result of the propagation of chaos linking the $N$-replicator dynamics and the McKean-Vlasov system (Section~\ref{sec:propchaos}). Then, for the $N$-replicator systems we investigate conditions that guarantee \textit{the persistence of the whole set of replicators} (Section~\ref{sec:persistence}). For the non-linear replicator,  we prove the existence and uniqueness of invariant measures for a particular type of interaction functions and noise (Section~\ref{section:invariant}). Section~\ref{sec:simulations} shows some simulations and numerical results related to the model treated in Section~\ref{section:invariant}. Finally, Section \ref{sec:conclusions} discuss our results and points out open questions regarding our research.


\section{Main results} \label{sec:framework}

Here we present our main definitions and assumptions and summarise the main results. 

\subsection{Notation and definitions}

\textbf{Notation}

\begin{itemize}
    \item $\mathbb{R}^d_{+}$ and $\mathbb{R}^d_{++}$ stand for the cones of real $d$-dimensional vectors with positive and strictly positive entries, respectively. A generic element of $\mathbb{R}^d_{+}$ or $\mathbb{R}^d_{++}$ will be denoted by bold letters $\xx, \yy, \zz$, etc., and when we need to stress its components, we will write $\xx= (x_1, x_2, \ldots,x_d)^\top$.
    \item For a metric space $(E, \rho)$ the notation $\mathcal{B}(E)$ stands for both, its Borel $\sigma$-field and the algebra of real-valued Borel measurable functions. Also, we write $\mathcal{P}(E)$ for the set of probability measures on $(E, \mathcal{B}(E))$.
    \item If $A$ is a subset of a topological space we write $\interior(A)$, $\overline A$, and $\partial A$ to denote its interior, closure, and boundary, respectively. In particular, $\interior (\simp{d})= \{\xx \in \simp{d}: \prod_{i=1}^d x_i \neq 0\}$ and $\partial \simp{d}= \{\xx \in \simp{d}: \prod_{i=1}^d x_i=0\}$.
    \item If $A, B$ are subsets of a topological space, the notation $A \Subset B$ means that $A$ is  compactly contained in $B$: $\overline A \subset B$, and $\overline A$ is compact.  
    \item For a function $f:E \to E'$, with $(E, d)$ and $(E', d')$ two metric spaces, the notation $\Lip (f)$ stands for the best Lipschitz constant of $f$, namely: $\Lip (f) \le K$ if and only if $d'(f(x), f(y)) \le K d (x, y)$ for every $x, y \in E$.
    \item If $E$ is a (regular enough) subset of $\mathbb{R}^d$,  we denote by $\mathcal{C} (E)$ the space of continuous real-valued functions on $E$. Analogously, $\mathcal{C}^k(E)$ stands for the $k$-times continuously differentiable real-valued functions on $E$. 
    \item If $E$ is a compact metric space and ${\bf A}:E \mapsto \mathbb{R}^d \otimes \mathbb{R}^{m}$ is a continuous map, we define the (uniform) Frobenius norm through:
\begin{align*}
||| {\bf A } |||  := \sup_{\xx \in E} \sqrt{\Tr{A^\top(\xx) A(\xx)} }. 
\end{align*}
    \item For a Banach space $E$ we denote with $\inner{\cdot}{\cdot}$ the duality product of $E^{*}$ with $E$. In particular, if $E= \mathbb{R}^d$, the notation $\inner{\xx}{\yy}$ stands for the usual inner product of vector of $\mathbb{R}^d$. 
    \item For $f \in \mathcal{B}_b(E)$ and $\mu$ a $\sigma$-finite measure on $E$,  we write $\mu f$ to denote the integral $\int_{E} f (\xx) \mu (d\xx)$. Regarding the point above, we also write $\inner{\mu}{f}$ with the same meaning.
    \item For probability measures $\mu, \nu$ on a measurable space $(E, \mathcal{E})$, the notation $d_{TV}(\mu,\nu)$ stands for total variation distance between $\mu$ and $\nu$, namely:
    \begin{align*}
     d_{TV} (\mu, \nu) = \sup_{A \in \mathcal{E}} \vert \mu (A)- \nu (A)\vert.
    \end{align*}
    If $E= \R{d}$, then $\Was_p (\mu, \nu)$ stands for the Wasserstein distance with respect to the underlying $p$-th vector norm, namely:
    \begin{align*}
    \Was_p (\mu, \nu) = \inf_{ \pi \in \Couplings (\mu, \nu)} \left ( \int_{\R{d} \times \R{d}} \Vert \xx - \yy \Vert^p \pi (d \xx, d\yy) \right )^{1/p}, 
    \end{align*}
    whenever the above quantity exists. Here, $\Couplings(\mu, \nu)$ is the family of probability measure on the product space $\R{d} \times\R{d}$ that have first and second marginals $\mu$ and $\nu$, respectively. 
    
\end{itemize}

\subsection{The stochastic replicator}

The structure of the SDEs (and later, McKean-Vlasov type equations) that we are going to consider take their general form from a geometric consideration. Recall that the Shahshahani metric on $\mathbb{R}^d_{++}$ is given by the metric tensor:
\begin{equation}
    s_{i,j}(\xx)=\frac{\norm{\xx}_1}{x_i}\delta_{i,j},\;(\xx\in E,\;i,j=1,\ldots,d).
\end{equation}
(see \cite{Shahshahani} for further details). The gradient of a $C^1$-function $U:\mathbb{R}^d\setminus\{{\bf 0}\}\to\mathbb{R}$ is transformed by the Shahshahani metric into $S(\xx)^{-1}\nabla U(\xx)$, where $S^{-1}(\xx)$ is the inverse of the diagonal matrix $S(\xx)$ whose entries are $s_{i,j}$, that is
\[S(\xx)^{-1}\nabla U(\xx)=\left(\begin{array}{c}
\frac{x_1}{\norm{\xx}_1}\frac{\partial}{\partial x_1}U(\xx)\\\vdots\\\frac{x_d}{\norm{\xx}_1}\frac{\partial}{\partial x_d}U(\xx)\end{array}\right)=\frac{1}{\norm{\xx}_1}\xx\circ \nabla U(\xx),\]
where $\circ$ denotes Hadamard product. In particular, for $\xx\in\simp{d}$ one gets $S(\xx)^{-1}\nabla U(\xx)=\xx\circ \nabla U(\xx)$. This gradient projected onto the tangent bundle at $\xx$ (hereafter, $\TT_\xx= \TT$) gives:
\begin{equation}\label{eq:grad_proj}
    \xx\circ\nabla U(\xx)-\langle \xx\circ\nabla U(\xx),\1\rangle \xx=\xx\circ(\nabla U(\xx)-\langle \xx,\nabla U(\xx)\rangle \mathbf{1}).
\end{equation}
Thus, at least if $A$ in equation \eqref{eq:original_replicator} is a symmetric payoff matrix, the original deterministic replicator has the form of a gradient-type dynamics on the simplex. 

Notice, however, that the right-hand side of equation \eqref{eq:grad_proj} makes sense if we replace $\nabla U$ by any $\mathbb{R}^d$-valued function $\bf F$. From now on, write:
\begin{align*}
\proj {\bf F} (\xx):= \xx \circ ({\bf F} (\xx)- \inner{{\bf F}(\xx)}{\xx}\mathbf{1}),
\end{align*}
and extend in an obvious way the definition for matrix-valued functions ${\bf G}: \mathbb{R}^d \mapsto \mathbb{R}^d \otimes \mathbb{R}^m$. 

Let $\Phi: \simp{d} \mapsto \mathbb{R}^d$ and $\Psi: \simp{d} \mapsto \mathbb{R}^{d} \otimes \mathbb{R}^d$ be given. We assume that, at least, both $\Phi$ and $\Psi$ are $\mathcal{C}^1 (\simp{d})$ functions. This allows us to affirm the existence of a uniquely defined matrix-valued function ${\bf A}_\Psi$ such that $\proj \Psi (\xx)= \xx \circ {\bf A}_\Psi (\xx)$ for every $\xx \in \simp{d}$. Consider the $d$-dimensional SDE given in \eqref{eq:eq00},
with suitable initial conditions on $\simp{d}$. We note first that this equation defines a $\simp{d}$-valued process that can be obtained via the projection on the simplex of a simple model of population dynamics (see the Appendix for the details of this derivation), and thus the paths remain in $\simp{d}$ by definition. Note also that by our assumptions, both drift and diffusion components are Lipschitz on $\simp{d}$, and thus it is well known that there exists a unique (pathwise) solution to the SDE \eqref{eq:eq00} which defines a well-posed $\simp{d}$-valued Markov process.  Let $(P^{(\Phi, \Psi)}_t: t \ge 0)$ be the semigroup associated to the Markov process $(\XX_t: t \ge 0)$ associated to the characteristics $(\Phi, \Psi)$. As usual, we consider the semigroup acting from the left on the space of bounded Borel function and acting from the right on Borel measures on $\simp{d}$. Let $B$ be a Borel set of $\simp{d}$. It is a routine fact that $(P^{(\Phi, \Psi)}_t: t \ge 0)$ is a Feller semigroup. By compactness, it admits invariant measures \cite[Theorem 12.39]{da2014introduction}. Recall that a probability measure $\mu \in \mathcal{P} (B)$ is invariant for the semigroup (equivalently, invariant for the process) if $\mu P^{(\Phi, \Psi)}_t = \mu$ for every $t \ge 0$, and in this case we write $\mu \in \mathcal{P}_{inv}(\Phi, \Psi; B)$. Let $\mathcal{P}_{erg}(\Phi, \Psi; B)$ be the subset of $\mathcal{P}_{inv} (\Phi, \Psi; B)$ of extremal invariant probability measures supported on $B$ (that is, ergodic probability measures that put no mass on $\simp{d} \setminus B$).

It is well known that for well-behaved Markov processes (for example, for processes whose semigroups enjoy the strong Feller property), different ergodic measures have disjoint supports. Thus, in order to obtain the uniqueness of an ergodic probability measure supported on the interior of the simplex, it seems reasonable to impose conditions on the noise in such a way that no submanifold of dimension less than $d-1$ is kept invariant by the dynamics. So, we will impose the following non-degeneracy of the noise.
\begin{assumptions} \label{ass:non_degenerate_noise}
\begin{align}\label{eq:non_degeneracy}
	\textbf{span} (\proj \Psi (\xx) _{\cdot, j}, j=1, \ldots, m)= \TT_\xx,
\end{align} 
for all $\xx\in \interior(\simp{d})$.
\end{assumptions}

Following \cite{benaim2008robust}, to assess the stochastic persistence of the process $\XX$ (and later, the stochastic persistence of interacting replicators), we are going to impose conditions on the natural deterministic process associated with it (namely, the process with no diffusion term). Define the \textbf{invasion growth rate}:

\begin{align} \label{eq:inv_rates}
h_{i} (\xx)= \Phi (\xx)_i - \langle \xx, \Phi(\xx) \rangle,
\end{align}
and set ${\bf h}:= (h_1, h_2, \ldots, h_d)^\top$.  Observe that the $i$-th invasion rate is just the exponent associated to the $i$-th population in the deterministic dynamics. Now, if $\mu \in \mathcal{P}_{inv} (\Phi, 0; \simp{d})$ is an invariant probability measure for the deterministic replicator under which type $i$ is absent, then $\mu$ can be regarded as an equilibrium distribution for a given strict population with strategies or types indexed by $I \subset \{1, 2, \ldots, d\}$ with $i \notin I$. Consider the following condition.

\bigskip

For every $\mu \in \mathcal{P}_{erg} (\Phi, 0; \partial \simp{d})$: 
\begin{align*}
    \max_{i=1, \ldots, d} \mu h_i > 0.
\end{align*}

The biological interpretation of this condition is just the requirement that any strict population of types in ecological equilibrium (for the deterministic dynamics) allows for the invasion of at least one of the absent types. It turns out that the above condition is equivalent to the following assumption (see e.g. \cite{schreiber2011}). 
\begin{assumptions} \label{ass:positive_Lyapunov}
There exists ${\bf p}= (p_1, p_2, \ldots, p_d)^\top$ a vector of strictly positive components such that:
\begin{align*}
\inf_{\mu \in \mathcal{P}_{erg} (\Phi, 0; \partial \simp{d})} \mu \langle {\bf p}, {\bf h} \rangle =: \rho > 0.    
\end{align*}
\end{assumptions}

The next theorem is just Theorem 1 in \cite{benaim2008robust}, which we state here for future comparison. 

\begin{theorem}\label{th:persistence}
    Assume that  $\Phi$ is $\mathcal{C}^2(\simp{d})$, and Assumptions \ref{ass:non_degenerate_noise} and \ref{ass:positive_Lyapunov} hold. Then, for every $r > 0$, there exists $\varepsilon= \varepsilon (r) > 0$ such that for every $\Psi$ of class $\mathcal{C}^3 (\simp{d})$ such that $ |||{\bf A}_\Psi|||<\varepsilon$ the following holds:
    \begin{enumerate}
    \item $\mathcal{P}_{inv} (int (\simp{d}))$ is a singleton, i.e., there exists a unique invariant probability measure $\mu$ supported on the interior of the simplex for the system \eqref{eq:eq00}. Moreover $\mu$ is absolutely-continuous w.r.t. Lebesgue measure, and the function $\xx \mapsto dist (\xx, \partial \simp{d})^{-r}$ is $\mu$-integrable.
    
    \item There exist positive constants $C, \alpha$ such that:
    \begin{align*}
        d_{TV} (\delta_{\xx}P_t, \mu) \le \dfrac{C e^{-\alpha t}}{\dist (\xx, \partial \simp{d})^{r}}, \quad \xx \in int(\simp{d}).
    \end{align*}
    \end{enumerate}
\end{theorem}

In this article, we extend the previous result to the setup of a family of replicators interacting in a mean-field regime. The next subsection provides the details. 

\subsection{Mean-field interacting replicators and McKean-Vlasov evolutionary games} \label{sec:SDEs}

As a first step, we extend the definition of our parameters to account for the interaction among replicators. Let be given a bounded measurable map $\Upsilon: \simp{d} \times \simp{d} \mapsto \mathbb{R}^d$. Observe that:
\begin{align*}
    \proj \Upsilon (\xx, \mu) =  \mu \proj \Upsilon (\xx, \cdot) = \proj \mu \Upsilon (\xx, \cdot).
\end{align*}
 We assume the following regularity.
\begin{assumptions} \label{ass:regularity_coeff}
$\Upsilon$ is a $\mathcal{C}^1$ map.
\end{assumptions}

In what follows, sometimes we will need to refer to the joint effect of both the endogenous and the interaction functions. To that aim, we define:
\begin{align*}
{\bf b}:= \proj \Phi + \proj \Upsilon,
\end{align*}
which is  continuous as a map from $\simp{d} \times \mathcal{P} (\simp{d})$, where the relevant topology in $\mathcal{P} (\simp{d})$ is the weak topology of probability measures. With this notation, equation~\eqref{eq:MFSDE} becomes
\begin{equation}\label{eq:particlesys}
d\XX^{(N; i)}_t =  {\bf b} (\XX_t^{(N; i)}, \mu_t^{(N)})  dt + \proj \Psi (\XX_t^{(N; i)}) d \WW_t^{(i)}, \quad i=1, \ldots N, 
\end{equation}
and McKean-Vlasov equation~\eqref{eq:mckean_vlasov_0} becomes
\begin{align}\label{eq:mckean_vlasov}
    d\Xi_t= {\bf b} (\Xi_t, \mu_t)dt + \proj \Psi (\Xi_t)d\WW_t,
\end{align}
where $\mu_t:= \Law (\Xi_t)$. 

\begin{remark}
Our hypotheses entail that $\bf b$, when seen as a map from $\simp{d} \times \mathcal{P}(\simp{d})$, is Lipschitz continuous for any Wasserstein distance in the measure variable, in the sense that for every $p\ge 1$ there exists a constant $C(p, d)$ such that for every $\xx, \yy \in \simp{d}$ and $\mu, \nu \in \mathcal{P}(\simp{d})$:
\begin{align*}
\Vert {\bf b} (\xx, \mu)- {\bf b} (\yy, \nu)\Vert_p \le C(p,d)\bigg( \Vert \xx-\yy\Vert_{p}+ \Was_p (\mu, \nu) \bigg ). 
\end{align*}
\end{remark}

The following result is proved in the next section. 

\begin{proposition} 
Under the above assumptions, Equation \eqref{eq:mckean_vlasov} has a unique pathwise (hence strong) solution.   
\end{proposition}


Now we are in a position to settle the main results of this paper. We start with a proposition that links the interacting particle system and the McKean-Vlasov replicator.

\begin{theorem} \label{th:prop_chaos}
    Assume $\Phi$ is $\mathcal{C}^2(\simp{d})$, and Assumptions  \ref{ass:non_degenerate_noise} and \ref{ass:regularity_coeff} hold. Let $\tilde \mu$ be a probability measure on the simplex, and consider the $N$-particle system \eqref{eq:particlesys} started from $\tilde \mu^{\otimes N}$. Consider also the McKean-Vlasov replicator \eqref{eq:mckean_vlasov} started from $\tilde \mu$. Then the propagation-of-chaos property holds in the following sense: for any positive time $t$, the sequence of random probability measures $(\mu_{[0, t]}^{(N)})_{N \ge 1}$, converges almost surely in the topology of the weak convergence of probability measures, to the deterministic probability measure $\mu_{[0, t]}= \Law (\Xi_{[0, t]})$.  
\end{theorem}

It is well known that the thesis of the above theorem is equivalent to the fact that, in the limit $N\to \infty$, a sort of asymptotic independence of the particles
holds (see Proposition 2.2, (i), at \cite{sznitman86}). More precisely, under the above conditions, for any $T > 0$ and any finite set of functions $g_k\in \mathcal{C}_b (C ([0, T]; \simp{d}))$, with $\E^{(N)}$ the expectation operator associated to the $N$-replicator system, and $\E^{(MV)}$ the corresponding expectation associated to the McKean-Vlasov replicator, we have:
\begin{align*}
    \lim_{N \to \infty} \mathbb{E}^{(N)}_{\tilde \mu^{\otimes N}} \left ( \prod_{k=1}^m g_k (\XX^{(N; k)}_{[0, t]}) \right ) = \prod_{k=1}^m \mathbb{E}^{(MV)}_{\tilde \mu} (g_{k} (\Xi_{[0, t]})).
\end{align*}

Let $P^{(N)}$ be the semigroup associated with the $N$-particle system. Let $\mathcal{P}^N_{inv}$ (resp. $\mathcal{P}^N_{erg})$ be the set of invariant (resp. ergodic) probability measures for $P^{(N)}$. The next theorem, which extends Theorem \ref{th:persistence}, is proved in Section \ref{sec:persistence}. It can be paraphrased as: if the typical replicator of the assemblage satisfies Assumption \ref{ass:positive_Lyapunov}, then the whole community is persistent. 

\begin{theorem}\label{th:persistence_interacting}
        Assume that $\Phi$ is $\mathcal{C}^2(\simp{d})$, and Assumptions \ref{ass:non_degenerate_noise}, \ref{ass:positive_Lyapunov}, \ref{ass:regularity_coeff} hold. Fix $N \ge 1$, $r > 0$. Then there exists $ \varepsilon:= \varepsilon(r)>0$ such that for any $\Psi$ of class $\mathcal{C}^2 (\simp{d})$ such that $\Vert \Upsilon \Vert + |||{\bf A}_\Psi|||<\varepsilon$, the following holds:
    \begin{enumerate}
    \item $\mathcal{P}^{N}_{inv} (\interior (\simp{d})^N)$ is a singleton, i.e., the $N$-replicator systems admits a unique invariant probability measure $\mu$ supported on the interior of $(\simp{d})^N$.  Moreover $\mu$ is absolutely-continuous w.r.t. Lebesgue measure, and the function $\xx \mapsto \dist (\xx, \partial (\simp{d})^N)^{-r/N}$ is $\mu$-integrable.
    
    \item There exist positive constants $C, \alpha$ such that:
    \begin{align*}
        d_{TV} (\delta_{\xx}P^{(N)}_t, \mu) \le \dfrac{C e^{-\alpha t}}{\dist (\xx, \partial (\simp{d})^N)^{r/N}}, \quad \xx \in \interior((\simp{d})^N).
    \end{align*}
    \end{enumerate}
\end{theorem}

Finally, let us consider the case of $\Phi(\xx)= A\xx$ and $\Psi(\xx)= \sigma I$, and assume that for the constant $\sigma$ and the payoff matrix $A$ satisfy the conditions:
\begin{enumerate}
\item[\textbf{[C1]}] For every $i \neq j$:
\begin{align*}
a_{ij}+a_{ji}-a_{ii}-a_{jj}= \sigma^2.
\end{align*}
\item[\textbf{[C2]}] For the matrix $A'$ with entries $a'_{ij}=a_{ij}- \dfrac{1}{2} \sigma^2$, there exists a vector $\alpha=(\alpha_1, \ldots, \alpha_d)$ of strictly positive entries summing up to $1$ such that $A' \alpha$ is a multiple of $\mathbf{1}$. 
\end{enumerate}

In the case of small $\sigma$, condition \textbf{[C1]} implies the existence of fitness equivalence (the diagonal and off-diagonal fitness payoffs are similar to each other), in the sense that any perturbation of the fitness payoffs that results in that intraspecific and interspecific interactions compensate each other, will preserve the condition of neutrality. Notice that the condition of neutrality is assumed in the neutral theory of ecology (\cite{hubbell2001}) but here it is derived as a condition for system persistence. Further, a small $\sigma$ also ensures the coexistence of strategies, such that the process spends little time near the pure strategy state, and, as shown below,  the process admits a unique invariant distribution of Dirichlet type. Condition \textbf{[C2]} implies the existence of an interior Nash equilibrium.\\

Consider the particular case where $\Upsilon$ is given by:
\begin{align*}
\Upsilon (\xx, \mu) = \delta F (\mu) \xx,
\end{align*}
where $F: \mathcal{P}(\simp{d}) \mapsto \R{d}\otimes \R{d}$ is a Lipschitz function with respect to the $\Was_1$ distance that ranges in the space of skew-symmetric real matrices and $\delta \in \mathbb{R}$ is a parameter that represents the interaction strength of a replicator with the rest of the alike replicators in the community. Observe that, in this case, the condition \textbf{[C1]} is preserved, in the sense that if $A$ satisfies this condition, then for every positive time the payoff matrix associated with the dynamics of the McKean-Vlasov replicator with characteristics $(\Phi, \Upsilon, \Psi)$ also satisfies the condition. Thus, if we regard \textbf{[C1]} as a proxy for neutrality, the prescription for $\Upsilon$ can be seen as a definition that preserves neutrality in the extended set up of mean-field replicator entities. In Section \ref{section:invariant} we prove the existence of invariant distribution for a class of McKean-Vlasov replicators associated with these parameters.

\begin{theorem}\label{th:invariant}

Assume $A$ satisfies \textbf{[C1]} and \textbf{[C2]}. Then there exists an interaction strength $ \delta_0 > 0$ such that for every $\delta < \delta_0$, the stochastic McKean-Vlasov replicator with $\Phi(\xx)= \tilde A \xx$ and $\Psi(\xx)= \sigma \mathbf{I}_{d \times d}$ constant admits a unique invariant probability measure putting no mass on the boundary of the simplex. Moreover, this invariant measure is of the Dirichlet family.
\end{theorem}

Although the above result states the existence of unique invariant measures for the above class of McKean-Vlasov replicator, by now we have no rigorous arguments for the ergodicity of the system. So, in Section \ref{sec:simulations} we provide statistical evidence to support this fact (at least in the simple setup of the $\simp{2}$-Mckean-Vlasov replicator). The evidence is then confronted to a test that allows us to accept the hypothesis that, for $\delta > 0$ small enough, every initial condition is attracted to the $\text{Beta}(s, 1-s)$ distribution, with $s= \dfrac{a_{12}-a_{22}}{\sigma^2-\delta}$.

\section{Propagation of chaos for the mean-field replicators}\label{sec:propchaos}

In this section, we settle some auxiliary results regarding the well-posedness of the equations we are dealing with, and we prove the propagation of chaos for the $N$-replicator system (Theorem \ref{th:prop_chaos}).

\subsection{Elementary results}

We start with a proposition that regards the existence and uniqueness of solutions to the nonlinear McKean-Vlasov equation. Its proof is relatively standard, and consequently, it is deferred to the  Apprendix~\ref{app:uniqueness}.

\begin{proposition} \label{prop:pathwise_uniqueness}
Under Assumptions \ref{ass:non_degenerate_noise} and \ref{ass:regularity_coeff}, equation \eqref{eq:mckean_vlasov} has a unique pathwise (hence strong) solution.   
\end{proposition}

Observe that the process $\Xi$ is \emph{not} a Markov process. Indeed, the family of maps $(P_s: s \ge 0)$ acting on $\mathcal{B}_b (\simp{d})$ defined via 
\begin{align*}
P_s f (\xx) := \E_{\xx} (f (\Xi_s)), 
\end{align*}
do not determine a semigroup of linear operators, in neat contrast to the Markovian setup. However, the semigroup property still holds in the dual sense. Define the one-parameter family of operators: 
\begin{align} \label{eq:Qs}
    \mathcal{Q}_t: \mathcal{P} (\simp{d}) \mapsto \mathcal{P} (\simp{d})
\end{align}
given by the prescription: $\mathcal{Q}_t \mu = \Law_{\mu}(\Xi_t)$, the law of $\Xi_t$ when started from an initial condition distributed as $\mu$.  

\begin{proposition}\label{prop:semigroup_property}
The family of operators $(\mathcal{Q}_t: t \ge 0)$ enjoys the semigroup property:
\begin{align}
    \mathcal{Q}_{t+s}= \mathcal{Q}_t \mathcal{Q}_s,
\end{align}
for all $s,t\geq 0$.
\end{proposition}

\begin{remark}
Still, the operators $\mathcal{Q}_s$ are non-linear: for $\mu, \nu \in \mathcal{P}(\simp{d})$, $\lambda \in [0,1]$, in general:
\begin{align*}
(\lambda \mu + (1-\lambda) \nu)\mathcal{Q}_t \neq \lambda \mu\mathcal{Q}_t
+ (1-\lambda )\nu \mathcal{Q}_t.
\end{align*}
\end{remark}

\begin{pf}
Let $\mu \in \mathcal{P}(\simp{d})$, and $t, s \ge 0$. Put $\mu_t:= \mathcal{Q}_t \mu$. Consider random variables $\Xi_0$ and $\Xi_1$ distributed as $\mu$ and $\mu_t$ respectively, and the processes $(\Xi^{(\Xi_0)}_{t+s}: s \ge 0)$ and $(\Xi^{(\Xi_1)}_{s}: s\ge 0)$. By \eqref{eq:mckean_vlasov}, for $s \ge 0$, on the one hand we have: 
\begin{align*}
    \Xi_{t+s}^{(\mu)} & \eqLaw \Xi_{t+s}^{(\Xi_0)} = \Xi_0 + \int_{0}^t \tilde \E ({\bf b} ( \Xi^{(\Xi_0)}_r, \tilde \Xi^{(\Xi_0)}_r )) dr + \int_{0}^{t} \proj\Psi ( \Xi^{(\Xi_0)}_r)d \WW_r \\
    & \quad + \int_{0}^s  \tilde \E ({\bf b} ( \Xi^{(\Xi_0)}_{r+t}, \tilde \Xi^{(\Xi_0)}_{r+t}) )dr + \int_{0}^{s} \proj\Psi ( \Xi^{(\Xi_0)}_{t+r}) d \WW_{t+r} \\
    &\eqLaw \Xi_0 + \int_{0}^t\tilde \E ({\bf b} ( \Xi^{(\Xi_0)}_r, \tilde \Xi^{(\Xi_0)}_r ))  dr + \int_{0}^{t} \proj\Psi ( \Xi^{(\Xi_0)}_r)d \WW_r \\
    & \quad + \int_{0}^s \tilde \E ({\bf b} ( \Xi^{(\Xi_0)}_{r+t}, \tilde \Xi^{(\Xi_0)}_{r+t}) ) dr + \int_{0}^{s} \proj\Psi ( \Xi^{(\Xi_0)}_{t+r})d \tilde \WW_r.
\end{align*}
where $\tilde \WW$ is a Wiener process independent of $(\WW_{r}:0 \le r \le t)$ and $\tilde \E$ stands for expectation with respect to the law of an independent copy of the process $\Xi^{(\Xi_0)}$. Observe that the sum of the first three terms in the last equation is distributed as $\Xi_1$ and is independent of $\tilde \WW$. Writing: $Y_{\cdot}:=\Xi^{(\Xi_0)}_{t+\cdot}$, the last line becomes:
\begin{align*}
    & \eqLaw \Xi_1 + \int_{0}^s \tilde \E ({\bf b} ( Y_r, \tilde Y_r) )dr + \int_{0}^{s} \proj\Psi ( Y_r)d  \WW_r \\
    &= Y_s \\
    &\eqLaw \Xi_{s}^{(\mu_t)},
\end{align*}
where the first-to-last line follows from pathwise uniqueness of the solutions of \eqref{eq:mckean_vlasov} (Proposition \ref{prop:pathwise_uniqueness}), and the last line is obvious. Our claim follows.
\end{pf}

For future reference, we establish a simple result. For the interested reader, we present the proof in Appendix~\ref{app:lemma1}.

\begin{lemma}\label{lemma:non_vanishing}
For every $\Xi_0 \in \interior (\simp {d})$:
\begin{align*}
    \mathbb{P}^{(MV)}_{\Xi_0}\left( \min_{i=1, \ldots, d} \Xi_t^{(i)} > 0  \text{ for all } t \ge 0 \right )=1. 
\end{align*}
An analogous result holds for the $N$- particle system, namely: for every $\vec \xx_0 \in \interior ((\simp{d})^N)$,
\begin{align*}
    \mathbb{P}^{(N)}_{\vec \xx_0}\left( \min_{j=1, \ldots, N} \Vert \XX^{(N; j)}_t \Vert > 0  \text{ for all } t \ge 0 \right )=1. 
\end{align*}
In words: with probability $1$, in finite time there is no absorption on the boundary of $\simp{d}$ (for the McKean-Vlasov replicator) nor on the boundary of $(\simp{d})^{N}$ (for the $N$-particle system).
\end{lemma}

As a direct result of Lemma~\ref{lemma:non_vanishing}, we have that for any $\xx \in \interior (\simp{d})$ and $t > 0$
\begin{align*}
\mathcal{Q}_t\delta_\xx(\partial\, \simp{d}) = 0,
\end{align*}
i.e., the trajectories of the solution of~\eqref{eq:mckean_vlasov} do not hit the boundary of the simplex in finite time. The following result regards the existence of a density of the law of the nonlinear process.

\begin{proposition} \label{prop:absolute_continuity}
For any $\xx \in \interior (\simp{d})$ and $t > 0$, $\mathcal{Q}_t\delta_\xx \ll \text{Leb}_{d-1}$.
\end{proposition}
\begin{pf}
For $t \ge 0$ and $\yy \in \simp{d}$, define:
\begin{align*}
\tilde {\bf b} (\yy, t):= {\bf b}(\yy, \mathcal{Q}_t \delta_\xx),    
\end{align*}
and consider the following diffusion with time-dependent coefficients:
\begin{align*}
    \Xi_{t}^{(\xx, \delta_\xx)}= \xx + \int _{0}^t \tilde {\bf b}(\Xi_{s}^{(\xx, \delta_\xx)}, s) ds + \int_{0}^t \proj \Psi (\Xi_{s}^{(\xx, \delta_\xx})) d\WW_s.
\end{align*}
Observe that the coefficients are Lipschitz-continuous and bounded, and moreover, on account of Assumption \ref{ass:non_degenerate_noise}, for every compact $\mathcal{D} \Subset \interior(\simp{d})$:
\begin{align*}
    \inner{\zz}{\proj \Psi(\yy) \proj \Psi (\yy)^\top \zz} > \varepsilon \Vert \zz \Vert^2,
\end{align*}
for a certain constant $\varepsilon>0$, for every $\yy \in \mathcal{D}$ and $\zz \in \mathcal{T}$. Now, for any $\xx \in \interior (\simp{d})$ consider:
\begin{align*}
\tau_{\mathcal{D}}:= \inf\{t \ge 0: \Xi_t^{(\xx, \delta_\xx)} \notin \mathcal{D}, \xx \in \mathcal{D}\},
\end{align*}
and define the finite measure on $\simp{d}$:
\begin{align*}
m_t^{\mathcal{D}}(A) := \prob (\Xi_t^{(\xx, \delta_\xx)} \in A, t \le \tau_{\mathcal{D}}).
\end{align*}
By a classical piece of the theory of parabolic PDEs with uniformly elliptic coefficients (see for example Chapter 6 of \cite{friedman1983} together with the representation given by Theorem 3.46, pag. 207, of \cite{pardoux2014}), we deduce that the measure  $m_t^{\mathcal{D}}(\cdot)$ is absolutely-continuous with respect to the Lebesgue measure.

Consider next a sequence of compact subsets $\{\mathcal{D}_{\epsilon}\}_{\epsilon}$ in $\interior (\simp{d})$,  such that $\mathcal{D}_{\epsilon'}\subset\mathcal{D}_{\epsilon}$ if $\epsilon<\epsilon'$, and $\lim_{\epsilon\rightarrow0}\mathcal{D}_{\epsilon}=\interior (\simp{d})$. Analogously as we have done with $\mathcal{D}$ above, define $\tau_{\mathcal{D}_{\epsilon}}$ and $m_t^{\mathcal{D}_{\epsilon}}(\cdot)$.  Let $A$ be a measurable subset of $\interior(\simp{d})$ with $\text{Leb}_{d-1} (A)=0$. Then $m_t^{\mathcal{D}_{\epsilon}}(A)=0$ for every $\epsilon > 0$. Since the trajectories of the solution of~\eqref{eq:mckean_vlasov} do not hit the boundary of $\simp{d}$ in finite time, we have that $m_t^{\mathcal{D}_{\epsilon}} (A) \uparrow \delta_\xx \mathcal{Q}_t (A)$ as $\epsilon \downarrow 0$. This proves the claim.
\end{pf}

\begin{remark}
The result still holds for the $N$-replicator systems, i.e., for every initial condition in $\interior ((\simp{d})^N)$ the law of $\vec \XX_t$ is absolutely continuous w.r.t. the $N$-fold product of the Lebesgue measure. 
\end{remark}

\subsection{Convergence of the empirical laws} 

Fix $t > 0$. Let $\mu^{(N)}$ be the random probability measure on $\mathcal{C}([0, t]; \simp{d})$ given by~\eqref{eq:random_pm}.
Now we will prove that $(\mu^{(N)})_N$ almost surely converges, in the weak sense, to the deterministic limit $\mu:=\Law (\Xi_{[0, t]})$. Recall that this means: if $m_N \in \mathcal{P} (\mathcal{P} (\mathcal{C}([0, t]; \simp{d})))$ is the law of the random measure $\mu^{(N)}$, we aim to prove that for every bounded continuous function $F \in \mathcal{C}_b (\mathcal{P} (\mathcal{C}([0, t]; \simp{d})))$, we have:
\begin{align*}
\int_{\mathcal{P} (\mathcal{C}([0, t]; \simp{d}))} m_N (d\nu) F(\nu) \rightarrow F (\mu),
\end{align*}
as $N$ goes to infinity. If we write $\E^{(N)}$ for the expectation operator associated with the $N$-replicator system, the above convergence is equivalent to asking that:
\begin{align*}
    \E^{(N)} ((\mu^{(N)} f - \mu f)^2) \to 0,
\end{align*}
as $N \to \infty$, for every $f \in \mathcal{C} (\mathcal{C}([0, t]; \simp{d})$. In turn, it can be proved that this will follow if we are able to prove that the family of probability measures $\E^{(N)} (\mu^{(N)}):=\Law (\XX_{[0, t]}^{(N; 1)})$ converges in $\mathcal{P} (\mathcal{C} ([0,t]; \simp{d})$ (see \cite{sznitman91}, pages 177 et seq.; see also \cite{meleard96}, pages 66 et seq.). So, this is what we will do now.



Assume that the $N$-particle system $(\XX^{(N; 1)}_{\cdot},...,\XX^{(N; N)}_{\cdot})$ in \eqref{eq:particlesys} is issued from independent and identically distributed conditions $(\XX^{(N; 1)}_{0},...,\XX^{(N; N)}_{0})$, and are driven by independent $d$-dimensional standard Brownian motions $(\WW^{(1)}_{\cdot},...,\WW^{(N)}_{\cdot})$. Notice that this implies that the law of $(\XX^{(N; 1)}_{\cdot},...,\XX^{(N; N)}_{\cdot})$ is symmetric.

\begin{theorem}
(Propagation of Chaos). Let  $(\Xi^{(i)}: i=1, \ldots, N)$ be a family of independent copies of the non-linear process solution of \eqref{eq:mckean_vlasov}, such that $\Xi^{( i)}_{0}=\XX^{(N; i)}_{0}$ for every $i=1, \ldots, N$ almost surely, and both $\XX^{(N; i)}_{\cdot}$ and $\Xi^{(i)}_{\cdot}$ are driven by indistinguishable Brownian motions. Then, under the above assumptions, for any finite time horizon $[0,S]$ there exists a positive constant $K=K(S,d)$ such that for every  $i=1, \ldots, N$:
\begin{align*}
\E\left(\sup_{t\in[0,S]}\left\Vert\XX^{(N; i)}_{t}- \Xi^{(i)}_{t}\right\Vert^2\right)\leq\frac{K}{N}.
\end{align*}
\end{theorem}
\begin{pf}
Let $\XX^{(N; i)}_{\cdot}$ and $\Xi^{(i)}_{\cdot}$ be the processes specified above, and let $S>0$. Then, for any $t\in[0,S]$ we have:
\begin{align*}
\left\Vert\XX^{(N; i)}_{t}- \Xi^{(i)}_{t}\right\Vert^2 &\leq  2S\int_{0}^{t}\left\Vert{\bf b} (\XX_s^{(N; i)}, \mu_s^{(N)})-{\bf b} (\Xi_s^{(i)}, \mu_s)\right\Vert^2 ds \\
&\quad + 2\left\Vert\int_{0}^{t}[\proj\Psi(\XX^{(N; i)}_{s})-\proj\Psi(\Xi^{(i)}_{s})]d\WW^{(i)}_{s}\right\Vert^2.
\end{align*}
Now, for the first term of the right-hand side, by the Lipschitz condition  of ${\bf b}$, there exists a $K_1=K_1(S,d)>0$ such that:
\begin{align*}
&\int_{0}^{t}\left\Vert{\bf b} (\XX_s^{(N; i)}, \mu_s^{(N)})-{\bf b} (\Xi_s^{(i)}, \mu_s)\right\Vert^2ds \leq K_1\int_{0}^{t}\left\Vert\XX^{(N; i)}_{s}-\Xi^{(i)}_{s}\right\Vert^2ds\\  
&\qquad \qquad \qquad +K_1\int_{0}^{t}\left\Vert\frac{1}{N}\sum_{j=1}^{N}\XX^{(N; j)}_{s}-\int_{\simp{d}}x\mu_s(dx)\right\Vert^2ds.
\end{align*}
On the other hand, notice that the term $\int_{0}^{\cdot}[\proj\Psi(\XX^{(N; i)}_{s})-\proj\Psi(\Xi^{(i)}_{s})]d\WW^{(i)}_{s}$ is an $(\FF_t: t \ge 0)$ square-integrable martingale, and thus by Doob's Maximal inequality we have that
\begin{align*}
&\mathbb{E}\left(\sup_{t\in[0,S]}\left\Vert\int_{0}^{t}[\proj\Psi(\XX^{(N; i)}_{s})-\proj\Psi(\Xi^{(i)}_{s})]d\WW^{(i)}_{s}\right\Vert^2\right) \\
&\leq4\mathbb{E}\left(\left\Vert\int_{0}^{S}[\proj\Psi(\XX^{(N; i)}_{s})-\proj\Psi(\Xi^{(i)}_{s})]d\WW^{(i)}_{s}\right\Vert^2\right).
\end{align*}
By It\^o's isometry and the Lipschitz assumption on $\proj\Psi(\xx)$, there exists a constant $K_2=K_2(S,d)>0$ such that:
\begin{align*}
&\mathbb{E}\left(\left\Vert\int_{0}^{S}[\proj\Psi(\XX^{(N; i)}_{s})-\proj\Psi(\Xi^{(i)}_{s})]d\WW^{(i)}_{s}\right\Vert^2\right) \\
&\leq K_2\mathbb{E}\left(\int_{0}^{S}\left\Vert\XX^{(N; i)}_{s}-\Xi^{(i)}_{s}\right\Vert^2ds\right).
\end{align*}
Therefore, we have that: 
\begin{align*}
&\E\left(\sup_{t\in[0,S]}\left\Vert\XX^{(N; i)}_{t}- \Xi^{(i)}_{t}\right\Vert^2\right)\leq K_3\int_{0}^{S}\E\left(\sup_{u\leq s}\left\Vert\Vert\XX^{(N; i)}_{u}- \Xi^{(i)}_{u}\right\Vert^2\right)ds\\
& \qquad \qquad \qquad +K_3\int_{0}^{t}\E\left(\sup_{u\leq s}\left\Vert\frac{1}{N}\sum_{j=1}^{N}\XX^{(N; j)}_{u}-\int_{\simp{d}}x\mu_u(dx)\right\Vert^2\right)ds,
\end{align*}
for some a constant $K_3=K_3(S,K_1,K_2)$.

Recall that we have assumed that $(\Xi^{(1)}_{\cdot}, \ldots, \Xi^{(N)}_{\cdot})$ are $N$ i.i.d. copies obeying \eqref{eq:mckean_vlasov} such that for all $j=1,...,N$, $\Xi^{(j)}_{0}=\XX^{(N; j)}_{0}$ almost surely, and $\Xi^{(j)}_{\cdot}$ and $\XX^{(N; j)}_{\cdot}$ are driven by indistinguishable Brownian motions. Then:
\begin{align*}
&\int_{0}^{t}\E\left(\sup_{u\leq s}\left\Vert\frac{1}{N}\sum_{j=1}^{N}\XX^{(N; j)}_{u}-\int_{\simp{d}}x\mu_u(dx)\right\Vert^2\right)ds \\
&\leq\dfrac{2}{N^2}\int_{0}^{S}\E\left(\sup_{u\leq s}\left\Vert\sum_{j=1}^{N}(\XX^{(N; j)}_{u}-\Xi^{(j)}_{u})\right\Vert^2\right)ds \\
&\quad +\dfrac{2}{N^2}\int_{0}^{S}\E\left(\sup_{u\leq s}\left\Vert\sum_{j=1}^{N}(\Xi^{(j)}_{u}-\int_{\simp{d}}x\mu_u(dx))\right\Vert^2\right)ds.
\end{align*}
The last term above is $O(N^{-1})$ by the Law of Large Numbers and, on the other hand, we have that:
\begin{align*}
\dfrac{1}{N^2}\int_{0}^{S}\E\left(\sup_{u\leq s}\left\Vert\sum_{j=1}^{N}(\XX^{(N; j)}_{u}-\Xi^{(j)}_{u})\right\Vert^2\right)ds\leq\int_{0}^{S}\E\left(\sup_{u\leq s}\left\Vert\XX^{(N; i)}_{u}-\Xi^{(i)}_{u}\right\Vert^2\right)ds,
\end{align*}
since, by symmetry, the laws of the $(\XX^{(N; j)}_{\cdot}-\Xi^{(j)}_{\cdot})$'s are exchangeable.
Finally, by  applying Gr\"onwall's inequality, we obtain:
\begin{align*}
\E\left(\sup_{t\in[0,S]}\left\Vert\XX^{(N; i)}_{t}- \Xi^{(i)}_{t}\right\Vert^2\right)\leq O(1/N)\exp\{2K_3S\},
\end{align*}
and therefore the theorem follows.
\end{pf}

We deduce from the previous result that the sequence of empirical measures $(\mu^{(N)})_N$ converges to the deterministic limit $\mu:=\Law (\Xi_{[0, t]})$,  in the weak sense. This fact plus the property that the law of the $N$-replicator $(\XX^{(N; 1)}_{\cdot},...,\XX^{(N; N)}_{\cdot})$ is symmetric (or exchangeable), implies that, for any $g_1$, $g_2$, ...., $g_k$ continuous and bounded functions on $\simp{d}$, with $k\leq N$, we have:
\begin{align*}
\lim_{N\rightarrow\infty}\E^{(N)}_{\mu^{\otimes N}}\left(\prod_{i=1}^{k}g_k(\XX^{(N; \sigma(i))}_{\cdot})\right)=\prod_{i=1}^{k}\inner{\mu}{g_i},
\end{align*}
where $(\XX^{(N; \sigma(1))}_{\cdot},...,\XX^{(N; \sigma(N))}_{\cdot})$ corresponds to any permutation of the vector $(\XX^{(N; 1)}_{\cdot},...,\XX^{(N; N)}_{\cdot})$ (see Proposition 2.2.i) in \cite{sznitman91} or Proposition 4.2 in \cite{meleard96}). That is, the $N$-particle system in the limit $N\rightarrow\infty$ becomes a sequence of independent and identically distributed solutions of the McKean-Vlasov equation.

\section{Strong stochastic persistence for the N-replicator system} \label{sec:persistence}

In this section we address the problem of persistence of the Markov processes $(\vec \XX^{(N)}_t: t\ge 0)$. For $N \ge 1$, the notation $(\simp{d})^N$ stands for the $N$-fold product of $\simp{d}$. The space $(\simp{d})^N$ is endowed with the topology of the topological sum. Since $\simp{d}$ is Polish, the same holds for its $N$-fold sum, for the metric:
\begin{align*}
    \dist_N (\vec\xx, \vec\yy) := \max_i \dist (\xx^{(i)}, \yy^{(i)}), 
\end{align*}
where the notation $\dist$ in the right-hand side stands for the $L^1$ distance on the simplex. We further set:
\begin{align*}
\dist_N (\vec \xx, \partial (\simp{d})^N) := \min_{j} \dist (\xx^{(j)}, \partial \simp{d}).
\end{align*}

Let $P^{(N)}$ be the semigroup associated with the $N$-replicator system~\eqref{eq:MFSDE}. The action of the generator associated with the semigroup on smooth functions can be expressed as follows. Given $\xx^{(1)}, \ldots, \xx^{(N)}$ arbitrary vectors in $\mathbb{R}^d$, we will write $\vec{\xx}$ for the vector in $\mathbb{R}^{dN}$ resulting from the concatenation, in lexicographic order, of the components of the $\xx^{(i)}$'s.  For a given smooth function $\varphi: (\simp{d})^N \mapsto \mathbb{R}$, for $i=1, \ldots, N$ and $\vec \xx \in (\simp{d})^N$ we denote by $\varphi^{(i)}_{\vec \xx}$ the map $\simp{d} \ni \yy \mapsto \varphi (\xx^{(1)}, \xx^{(2)}, \ldots, \yy, \ldots, \xx^{(N)})$ where $\yy$ is in the place of the $i$-th block of order $d$. A straightforward computation shows that whenever $\varphi$ is a twice-continuously differentiable real-valued function, the action of the generator is given by:

\begin{align*}
 L^{(N)} \varphi (\vec \xx) &= \sum_{i=1}^N \inner{\nabla \varphi^{(i)}_{\vec \xx} (\xx^{(i)}) }{ \proj \Phi (\xx^{(i)})} \\
 & \qquad +\sum_{i=1}^N \inner{\nabla \varphi^{(i)}_{\vec \xx} (\xx^{(i)}) }{ \proj \Upsilon (\xx^{(i)}, \mu_{\vec \xx})} \\
    & \qquad + \dfrac{1}{2} \sum_{i=1}^N \Tr  \left (\proj \Psi(\xx^{(i)})^\top \nabla \nabla ^\top \varphi^{(i)}_{\vec \xx}(\xx^{(i)}) \proj \Psi (\xx^{(i)}) \right ), 
\end{align*}
where: $\mu_{\vec{\xx}}= \dfrac{1}{N} \sum_{j=1}^N \delta_{\xx^{(i)}}$, and $\nabla \nabla^\top \varphi^{(i)}_{\vec \xx}$ is the Hessian matrix associated to $\varphi^{(i)}_{\vec \xx}$.

\begin{definition}\label{def:SSP} (see \cite{hening2018}, \cite{videla2022})
	We say that a $(\simp{d})^N$-valued, right-continuous  Markov process $(Y_{t}: t \ge 0)$ is {\bf strongly stochastically persistent} (SSP) if there exists a unique invariant probability measure $\pi$ such that $\pi (\boundary((\simp{d})^N))=0$ and such that for every $\xx \in \interior((\simp{d})^N)$:
	\begin{align*}
	d_{TV}(\prob_{x} (Y_t \in \cdot), \pi (\cdot)) \rightarrow 0 , 
	\end{align*}   
	as $t \rightarrow \infty$.
\end{definition}

Here we will prove that under our running assumptions, $(\vec \XX^{(N)}_t: t\ge 0)$ is a SSP processes taking values on $(\simp{d})^N$. In fact, we will prove an extension of Theorem \ref{th:persistence} and, in particular, we will obtain geometric rates of convergence for both processes.

To this end, we introduce some notation. We will say that $\xx \in \interior(\simp{d})$ is $\varepsilon$-close to extinction if there exists $i \in I$ such that $x_i < \varepsilon$. For given $\varepsilon$, define:
\begin{align*}
\text{Ext}(\varepsilon):=\{\xx \text { is } \varepsilon -\text{close to extinction}\}.
\end{align*}
We write $\text{Ext}_N(\varepsilon):= \{\vec \xx \in (\simp{d})^N: \xx^{(i)} \in \text{Ext}(\varepsilon) \text{ for some } i\}$. 

To prove convergence to an invariant measure for $(\vec \XX^{(N)}_t:t \ge 0)$, we first observe that it is enough to prove that for every $t> 0$ there is convergence for the chain $(\vec \XX^{(N)}_n:= \vec \XX^{(N)}_{tn}: n \ge 0)$, i.e. the chain sampled from $\vec \XX^{(N)}$ every $t$ time units. This claim follows from the Feller property of the semigroup $(P^{(N)}_t: t \ge 0)$ (which can be derived in a rather standard way from the Lipschitz continuity of the parameters) and from the following general fact about Feller processes.
 
\begin{lemma} \label{lemma:feller_uniqueness}
Let $(P_t: t \ge 0)$ be the Feller semigroup associated with a  Markov process $(Y_t: t \ge 0)$ taking values on a compact metric space $(E, \rho)$. Denote by $(P^{*}_t: t \ge 0)$ the dual semigroup acting on finite measures. Fix $F \in \mathcal{B}(E)$. Assume that for every $s \in (0,1]$ there exists a unique probability measure $\pi_s$ on $(E, \mathcal{B}(E))$ with $\pi_s (F^C)=0$ and such that $P^{*}_s\pi_s= \pi_s$. Then the process $Y$ has a unique invariant probability measure $\pi$ with $\pi(F^C)=0$.      
\end{lemma} 
 
\begin{pf}
For every $t \in (0,1] \cap \mathbb{Q}$, the unique invariant measure for $P^{*}_t$ is also invariant for $P^{*}_1$. Hence, $\pi_t= \pi_1$ for every rational $t$. Then, for general $t$ and any $n \in \mathbb{N}$, we have:
\begin{align*}
\pi_1 & = P^{*}_{\frac{\lfloor nt \rfloor}{n}} \pi_1,
\end{align*}
Since $P$ is Feller and $E$ is compact, $P^{*}$ is a weakly-continuous semigroup. Hence, sending $n$ to $\infty$ in the last equation, we obtain $\pi_1 f = P^{*}_t\pi_1 f$ for every $f \in \mathcal{C} (E)$. Since the class of continuous functions is measure-determining on $E$, we get that $\pi_1 = P^{*}_t \pi_1 $, and by the uniqueness assumption, $\pi_t=\pi_1$.
\end{pf}
Fix $t>0$, and let $P:=P^{(N)}_t$. Consider the following statement.

\bigskip

\begin{minipage}[t]{0.9\textwidth}
\textbf{[H]} There exists a continuous function $\mathcal{H}: \interior( (\simp{d})^N) \mapsto \mathbb{R}_+$ and constants $\alpha \in (0,1)$ and $C > 0$ such that:
\begin{enumerate}
    \item $\lim_{\varepsilon \to 0} \inf_{\vec \xx \in \text{Ext}_N (\varepsilon)}\mathcal{H} (\xx) = \infty$.
    \item $P \mathcal{H} (\vec \xx) \le \alpha \mathcal{H} (\vec \xx) + C$ for every $\vec \xx \in \interior ((\simp{d})^N)$. 
\end{enumerate}
\end{minipage}

\bigskip
 
Observe that $\textbf{[H]}$ entails that for any $\xx \in \interior ((\simp{d})^N)$, any weak limit point, say $\nu$, of the sequence of probability $(\delta_{\vec \xx} P^n)_n$ satisfies $\nu (\partial( (\simp{d})^N))=0$. For suppose this is not the case, and given $\vec \xx \in \interior((\simp{d})^N)$, assume that a limit point $\nu$ of $(\delta_{\vec \xx}P^n)_n$ satisfies $\nu (\partial ((\simp{d})^N))> 2 \varepsilon$ for a given $\varepsilon > 0$. Apply $n$ times the second condition of $\textbf{[H]}$ to get:
\begin{align*}
    P^n \mathcal{H} (\vec \xx) \le \alpha^n \mathcal{H} (\vec \xx) + D,
\end{align*}
where $D:=\dfrac{C}{1-\alpha}$. Define $\eta: = \dfrac{\mathcal{H} (\vec \xx) + D}{\varepsilon}$ and consider the set $K_\varepsilon:= \mathcal{H}^{-1} ([0, \eta])$. Since $\mathcal{H}$ is continuous, $K_\varepsilon$ is a closed (indeed, compact) set contained in $\interior((\simp{d})^N)$ (by the first point at $\textbf{[H]}$). Then, for every $n$, Markov inequality yields:
\begin{align}\label{ineq:tightness}
    P^n \mathbf{1}_{K_\varepsilon^C} (\vec \xx) \le \varepsilon. 
\end{align}
Consider a continuous function $g: (\simp{d})^N \mapsto [0,1]$ with the property $g \equiv 1$ in $\partial ( (\simp{d})^N)$, $g\equiv 0$ on $K_\varepsilon$. Since $g \le \mathbf{1}_{K_\varepsilon^C}$, by \eqref{ineq:tightness} we obtain $\delta_{\vec \xx}P^n g  \le \varepsilon$ for every $n$. But $\nu g > 2 \varepsilon$ which contradicts the fact that $\nu$ is a weak-limit point. 

\begin{remark}
Consider for a while the dynamics given by the McKean-Vlasov replicator $\Xi$ discussed in the previous sections, and for every $t \ge 1$ consider the operator $\tilde P_t$ acting on bounded functions as $\tilde P_t f(\xx):= \E_{\xx} (f (\Xi_t))$. As remarked above, in general, we do not have $\tilde P_{nt} = \tilde P^{n}_t$. So, the previous argument breaks down from the very start in the case of McKean-Vlasov replicators, and thus new tools are needed to face the problem of persistence in this setup. 
\end{remark}

Thus, $\textbf{[H]}$ guarantees that the limit points of the sequence $(\Law(\vec \XX^{(N)}_{tn}): n \ge 0)$ support the whole ensemble of replicators. Also, observe that if $\nu$ is an invariant probability measure for the process $\vec \XX^{(N)}$ that supports the whole ensemble of replicators, then by definition $\nu P = \nu$. Let $A \subseteq \interior((\simp{d})^N)$ be a measurable set with $\text{Leb} (A)=0$. Then:
\begin{equation*}
    \nu (A)  =  \nu P (A) = \int_{(\simp{d})^N} \nu (d\vec \xx) \delta_{\vec \xx} P \mathbf{1}_{A} = \int_{\interior((\simp{d})^N)} \nu (d\vec \xx)\delta_{\vec \xx} P \mathbf{1}_{A},
\end{equation*}
since $\nu$ does not charge the boundary of $(\simp{d})^N$. By Proposition \ref{prop:absolute_continuity}, for every $\vec \xx$ there exists a density, say $p(\vec \xx, \cdot)$, of the probability measure $\delta_{\vec \xx} P$ with respect to the Lebesgue measure. Consequently:
\begin{multline*}
    \nu (A) = \int_{\interior((\simp{d})^N)} \nu (d\vec \xx)\delta_{\vec \xx} P \mathbf{1}_{A}
    \\= \int_{\interior(\simp{d})} \nu (d\vec \xx) \int_{A} p (\vec \xx, \vec \yy) d\vec \yy= \int_{A} d\vec \yy \left ( \int_{\interior((\simp{d}))^N } \nu (d\vec \xx) p (\vec \xx, \vec\yy) \right ), 
\end{multline*}
by Fubini's theorem. Since $\text{Leb} (A)=0$, the last integral vanishes, and we deduce $\nu (A)=0$. In other words, if $\nu$ is an invariant measure for the process $\vec  \XX^{(N)}$ under which the full system persists, then it is continuous with respect to the Lebesgue measure.  

Finally, assume that we are able to prove that there exists a unique invariant measure, say $\mu$, for the semigroup $P$ which is supported at the interior of $(\simp{d})^N$. We claim that $\mathcal{H}$, as in $\textbf{[H]}$, is $\mu$-integrable. Indeed, for $n \ge 0$, define $\mathcal{H}_n:= \mathcal{H} \wedge n$. Then, for every $n \ge 1$:
\begin{equation*}
    \mu \mathcal{H}_n = \lim_{k \to \infty} \mu P^k \mathcal{H}_n \le \mu (\limsup_{k \to \infty} P^k \mathcal{H}_n) \le D,  
\end{equation*}
where the equality follows since $\mu$ is invariant for $P$, and the first inequality is (reverse) Fatou's Lemma.  Apply (direct) Fatou's Lemma to obtain:
\begin{equation*}
    \mu \mathcal{H}  = \int_{\interior((\simp{d})^N)} \lim_{n \to \infty} \mathcal{H}_n (\vec \xx) \mu (d\xx)
     \le \liminf_{n \to \infty} \int_{\interior((\simp{d})^N)}  \mathcal{H}_n (\vec \xx) \mu (d\vec \xx)
    \le D,
\end{equation*}
and we deduce that $\mathcal{H}$ is $\mu$-integrable. 

Actually, the existence of a function $\mathcal{H}$ as in $\textbf{[H]}$ implies stronger results (see \cite{benaim2008robust}, page 186, and the references therein, for the following lemma). 

\begin{lemma} \label{lemma:convergence_lemma}
Assume that $\mathcal{H}$ exists as in $\textbf{[H]}$. Assume also that $P$ possesses a density with respect to the Lebesgue measure. Then, there exists a unique invariant probability measure $\pi$ for the kernel $P$, and moreover: 
\begin{enumerate}
\item $\pi \ll \text{Leb}$.
\item $\mathcal{H}$ is $\pi$-integrable.
\item There exists $\rho \in (0,1)$ and $C' > 0$ such that:
\begin{align}\label{eq:geometric_rate}
d_{TV} (\delta_{\vec \xx} P^n, \pi) \le C' (1+\mathcal{H}(\vec \xx))\rho ^n, \quad \vec \xx \in \interior((\simp{d})^N).
\end{align}
\end{enumerate}
\end{lemma}

\begin{remark}
Observe that \eqref{eq:geometric_rate} implies that $\pi$ attracts every initial measure $\nu$ with $\nu (\partial (\simp{d})^N)=0$. Indeed, for any compact set $K \subset \interior{(\simp{d})^N}$, by the definition of total variation norm and Fubini's theorem:
\begin{align*}
    \Vert \nu P^n - \pi \Vert_{TV} &= \sup_{f \textbf{ Borel }, \vert f \vert \le 1 } \left ( \int_{(\simp{d})^N} \nu (d\vec\xx)P^n f(\vec\xx)- \int_{(\simp{d})^N} \pi (d\vec\yy) f(\vec\yy)\right ) \\
    &= \sup_{f \textbf{ Borel }, \vert f \vert \le 1 }  \int_{(\simp{d})^N}  \nu (d\vec\xx) \int_{(\simp{d})^N} \left(P^n (\vec\xx, d\vec\yy) f(\vec\yy)- \pi (d\vec\yy) f(\vec\yy)\right) \\
    & \le  \int_{(\simp{d})^N}  \nu (d\vec\xx)  \sup_{f \textbf{ Borel }, \vert f \vert \le 1 }  \int_{(\simp{d})^N} \left(P^n (\vec\xx, d\vec\yy) f(\vec\yy)-   \pi (d\vec\yy) f(\vec\yy) \right) \\
    & \le \int_{K}  \nu (d\vec\xx) \Vert \delta_{\vec\xx} P^n -\pi \Vert_{TV} + 2 \nu (\interior((\simp{d})^N) \setminus K),
\end{align*}
since the total variation norm is bounded by $2$ and $\nu$ does not charge the boundary. Since $(\simp{d})^N$ is Polish, $\nu$ is regular, and then given $\varepsilon> 0$ we can choose $K_\varepsilon$ such that the last term is smaller than $\varepsilon$. The result follows since the term inside the integral converges uniformly to $0$ on $K_\varepsilon$ as $n \to \infty$.
\end{remark}

\begin{theorem} \label{th:persistence_3}
 Assume that  $\Phi$ is $\mathcal{C}^2(\simp{d})$, and Assumptions \ref{ass:non_degenerate_noise}-\ref{ass:positive_Lyapunov}-\ref{ass:regularity_coeff} hold. Then for every $r > 0$ there exists $\delta > 0$ such that for any $\Upsilon, \Psi$ of class $\mathcal{C}^2$  such that $\Vert\proj \Upsilon \Vert +|||\A_{\Psi}|||<\delta$, the following hold:
 \begin{enumerate}
 \item $\mathcal{P}_{inv}^{(N)} (\interior ((\simp{d})^N))= \{\nu\}$; i.e., there exists a unique invariant probability measure supported on the interior of the $N$-fold topological sum of $\simp{d}$ for the system \eqref{eq:particlesys}. Moreover $\nu$ is absolutely continuous w.r.t. Lebesgue measure, and the function $\interior((\simp{d})^N) \ni \vec \xx \mapsto \dist_N (\vec \xx, \partial ((\simp{d})^N))^{-r/N}$ is $\nu$-integrable.
 \item There exist positive constants $C', \alpha'$ such that:
    \begin{align*}
        \Vert \delta_{\xx}P^{(N)}_t- \nu \Vert_{TV} \le \dfrac{C'  e^{-\alpha' t}}{\dist_N (\xx,  \partial((\simp{d})^N))^{r/N}}, \quad \xx \in \interior ((\simp{d})^N).
    \end{align*}
\end{enumerate}
\end{theorem}

\begin{pf} (of Theorem \ref{th:persistence_3}).
We follow closely the steps in the proof Theorem 3.1 of \cite{benaim2008robust}. Recall the definition of the invasion rates \eqref{eq:inv_rates}. There exists $\alpha>0$, a neighbourhood $\mathcal{U}$ of $\partial \simp{d}$ and a function $W: \simp{d} \mapsto \mathbb{R}$ of class $\mathcal{C}^2 (\simp{d})$ such that:
\begin{align}\label{eq:exists_alpha}
    \inner{{\bf p}}{{\bf h}(\xx)} + \inner{\nabla W(\xx)}{\proj \Phi(\xx)} > \alpha,
\end{align}
for every $\xx \in \mathcal{U}$ (see \cite{Hofbauer2003b}, Theorems 3.4 and 4,4, and Remark 3.5), and where $\bf p$ is taken as in Assumption \ref{ass:positive_Lyapunov}.
Define the map $V_N:(\mathcal{U}\setminus\partial \simp{d}))^N\mapsto \mathbb{R}$  given by
\begin{align*}
    V_N(\vec \xx):= \dfrac{1}{N}\sum_{i=1}^N V (\xx^{(i)}),
\end{align*}
where $V(\xx):= \sum_{j=1}^d p_j \ln (x_j) + W (\xx)$. Since ${\bf h}(\xx)=\text{diag}(\xx)^{-1}\proj\Phi(\xx)$, we have: 
\begin{equation*}
\sum_{i=1}^N \inner{(\nabla V_N(\vec \xx))^{(i)}}{\proj \Phi (\xx^{(i)})}  = \dfrac{1}{N}  \sum_{i=1}^N \inner{\nabla V (\xx^{(i)})}{\proj \Phi (\xx^{(i)})} 
 > \alpha,
\end{equation*} 
on $(\mathcal{U}\setminus \partial\simp{d})^N$. Apply the generator $L^{(N)}$ to $V_N$ to obtain:
\begin{align*}
    L^{(N)} V_N(\vec \xx) & = \sum_{i=1}^N \inner{\nabla V^{(i)}_{N, \vec \xx}(\xx^{(i)})}{\proj \Phi (\xx^{(i)})} \\
    & \quad +\sum_{i=1}^N \inner{\nabla V^{(i)}_{N, \vec\xx} (\xx^{(i)}) }{ \proj \Upsilon (\xx^{(i)}, \mu_{\vec \xx})} \\
    & \qquad + \dfrac{1}{2} \sum_{i=1}^N \Tr  \left (\proj \Psi(\xx^{(i)})^\top \nabla \nabla ^\top V^{(i)}_{N, \vec \xx}(\xx^{(i)}) \proj \Psi (\xx^{(i)}) \right ),
\end{align*} 
where the map $V^{(i)}_{N,\vec \xx}$ is given by $\simp{d} \ni \yy \mapsto V_N(\xx^{(1)}, \xx^{(2)}, \ldots, \yy, \ldots, \xx^{(N)})$, where $\yy$ is in the place of the $i$-th block of order $d$. 

Since $W$ is of class $\mathcal{C}^2 (\simp{d})$, the above expression can be made greater than $\alpha/2$, say, for $\vec\xx \in (\mathcal{U}\setminus\partial \simp{d}))^N$ by choosing $\tilde \varepsilon_1 > 0$ small enough and imposing $\Vert \proj \Upsilon \Vert + ||| \A_{\Psi} ||| <\tilde \varepsilon_1$.

Define $\lambda:= \dfrac{r}{\inf p_i}$, and consider the function $\mathcal{H}: (\interior(\simp{d}))^N \mapsto \mathbb{R}_{++}$ given by $\mathcal{H}: = \exp (-\lambda V_N)$, which is of class $\mathcal{C}^2 ((\interior (\simp{d}))^N)$. Furthermore
\begin{align*}
\lim_{\eta \to 0} \inf_{\{\vec\zz:\zz^{(i)} \in \text{Ext} (\eta) \setminus \partial \simp{d}\text{, } \forall i=1,...,N\}} \mathcal{H} (\vec\zz) =  \infty.
\end{align*}

Since $W$ is bounded on $\simp{d}$, with $K:= \exp (-\lambda \inf_{\zz \in \simp{d}} W (\zz))$, we have:
\begin{align} \label{eq:integrability_at_limit}
 \mathcal{H} (\vec\xx) &\ge  K \exp (-\lambda(1/N)\sum_{i=1}^{N} \sum_{j=1}^{d} p_j\ln (\xx^{(i)}_j)) \nonumber \\
 & \ge K\prod_{i=1}^{N} \exp (-(r/N) \sum_{j=1}^{d}\ln (\xx^{(i)}_j))\nonumber \\
 &= K\prod_{i=1}^{N}\dfrac{1}{(\xx^{(i)}_1 \xx^{(i)}_2\ldots \xx^{(i)}_d)^{r/N}} \nonumber \\
 &\ge K\prod_{i=1}^{N}\dfrac{1}{\dist(\xx^{(i)},\partial (\simp{d}))^{r/N}}\nonumber \\
 &\ge \dfrac{K}{\dist_N (\vec\xx, \partial (\simp{d})^N)^{r/N}}.
\end{align}
Let $\Gamma^{(N)}$ be the \emph{carr\'e du champ} operator associated to the operator $L^{(N)}$, whose action on smooth functions $\varphi: (\simp{d})^N \mapsto \mathbb{R}$ is given by: 
\begin{align*}
\Gamma^{(N)} \varphi (\vec \xx) = \dfrac{1}{2} \sum_{i=1}^N \sum_{j,k} [\proj \Psi (\xx^{(i)}) \proj \Psi (\xx^{(i)})^\top]_{j,k} \partial_{j} \varphi_{\vec \xx} (\xx^{(i)}) \partial_{k} \varphi_{ \vec \xx} (\xx^{(i)}).
\end{align*}
For smooth functions $f: \mathbb{R} \mapsto \mathbb{R}, \varphi: (\simp{d})^N \mapsto \mathbb{R}$, we have:
\begin{align*}
L^{(N)} f(\varphi (\vec \xx))= f' (\varphi (\vec \xx)) L^{(N)} \varphi (\vec \xx) + f'' (\varphi (\vec \xx)) \Gamma^{(N)} \varphi (\vec \xx). 
\end{align*}
Hence:
\begin{align*}
L^{(N)}\mathcal{H} (\vec \xx)= - \lambda \mathcal{H} (\vec \xx) \left ( L^{(N)}  V_N (\vec \xx) -\lambda \Gamma^{(N)} V_N (\vec \xx) \right),
\end{align*}
and then we can choose $\tilde \varepsilon < \tilde \varepsilon_1$ and impose $\Vert \proj \Upsilon\Vert+|||\A_\Psi|||<\tilde \varepsilon$  such that the term inside the parentheses is not smaller than $\alpha/3$, say, on $(\mathcal{U} \setminus \partial \simp{d})^N$. And then, for this choice:
\begin{align*}
    L^{(N)} H (\vec\xx) \le -\dfrac{\lambda\alpha}{3} H (\vec\xx),  
\end{align*}
on $(\mathcal{U} \setminus \partial \simp{d})^{N}$. Observe that $L^{(N)} \mathcal{H}$ is bounded on $(\interior (\simp{d}) \setminus \mathcal{U})^N$ by some absolute constant $C$. Thus, with $\beta := \dfrac{\lambda\alpha}{3}$, we obtain:
\begin{align}\label{eq:bound_generator}
    L^{(N)} \mathcal{H}(\vec\xx) \le -\beta \mathcal{H} (\vec\xx) + C,
\end{align}
for all $\vec\xx \in (\interior(\simp{d}))^N$.

Let  $\eta_n:= \inf\{t\geq 0:\mathcal{H} (\vec\XX^{(N)}_t) > n\}$, and observe that since $\mathcal{H}$ is bounded on compact sets of $(\interior(\simp{d}))^N$,  by Lemma \ref{lemma:non_vanishing}, we have $\eta_n \to \infty$ a.s. For given $t \ge 0$, set $t_n:= t \wedge \eta_n$. For every $n\ge 1$, the process:
\begin{align*}
M^{n}_t = e^{\beta t_n}\mathcal{H} (\XX^{(N)}_{t_n}) - \mathcal{H}(\vec\XX^{(N)}_{t_0}) - \int_{0}^{t_n} e^{\beta s} (\beta \mathcal{H} (\vec\XX^{(N)}_s)+ L^{(N)} \mathcal{H} (\vec\XX^{(N)}_s))ds, \quad t \ge 0,
\end{align*}
is a martingale, and thus for every $n\ge 1$, by \eqref{eq:bound_generator}:
\begin{equation*}
    \E_{\vec \xx} (e^{\beta t_n} \mathcal{H} (\vec\XX^{(N)}_{t_n}))= \mathcal{H} (\vec \xx) + \E \left ( \int_{0}^{t_n} C e^{\beta s}\right )
    \le \mathcal{H} (\vec \xx) + \dfrac{C}{\beta} (1- e^{\beta t}). 
\end{equation*}
Hence:
\begin{align*}
\E_{\vec \xx}(e^{\beta t} \mathcal{H} (\vec\XX^{(N)}_{t}) \mathbf{1}_{t \le t_n}) \le \mathcal{H} (\vec \xx) + \dfrac{C}{\beta}, 
\end{align*}
and then, 
\begin{align*}
    \E_{\vec \xx} (\mathcal{H} (\vec\XX^{(N)}_{t}) \mathbf{1}_{t \le t_n}) \le e^{-\beta t} \mathcal{H} (\vec \xx) + \dfrac{C}{\beta}.
\end{align*}
Since $t_n \to t$ a.s., by Monotone Convergence, we deduce that for every $t \ge 0$ and $\vec\xx \in (\interior(\simp{d}))^N$:
\begin{align*}
    \E_{\vec \xx} (\mathcal{H} (\vec\XX^{(N)}_{t})) \le e^{-\beta t} \mathcal{H} (\vec \xx) + D, 
\end{align*}
where $D$ is a constant. In particular, there exists $\gamma \in (0, 1)$ such that, at time $1$:
\begin{align*}
    P^{(N)}_1 \mathcal{H} \le \gamma \mathcal{H} + D, 
\end{align*}
on $(\interior (\simp{d}))^N$. By Proposition \ref{prop:absolute_continuity}, $P^{(N)}_1$ possesses a density with respect to the Lebesgue measure. Hence, Lemma \ref{lemma:convergence_lemma} applies to $P^{(N)}_1$ and we deduce the existence of a unique invariant probability measure $\mu$ for $P^{(N)}_1$ such that $\mu (\partial (\simp{d})^N)=0$, and moreover $\mu$ attracts every initial condition on $(\interior(\simp{d}))^N$ at geometric rate. Of course, the same line of reasoning applies to $P^{(N)}_t$ for every $t \in (0,1]$, and we deduce that for every $t \in (0,1]$ there exists a unique invariant measure, say $\pi_t$, for the chain $(\vec\XX^{(N)}_{nt}:n \ge 0)$. The uniqueness follows from Lemma \ref{lemma:feller_uniqueness}. 

The remaining parts of the result are exactly those in the reference \cite{benaim2008robust}.   
\end{pf}

\section{Invariant distributions for the McKean-Vlasov stochastic replicator}\label{section:invariant}

We go back to the initial model of Fudenberg and Harris, as treated in \cite{hofbauerimhof2009}, and consider a non-linear, McKean-Vlasov-type extension of the model. Assume first that $d=2$, that $A$ is given, and that $\Sigma=\sigma I$. Recall that in this case, each replicator solves~\eqref{eq:stochastic_replicator_1}. From Theorem 3.6 of \cite{hofbauerimhof2009} we deduce the following result.

\begin{proposition}
Assume that $a_{12}-a_{22}> 0$ and $a_{21}-a_{11}> 0$, and furthermore:
\begin{align*}
a_{12}+a_{21}-a_{11}-a_{22}= \sigma^2.
\end{align*}
Then, the system \eqref{eq:stochastic_replicator_1} has a unique invariant distribution that does not charge the boundary, namely the absolutely continuous measure with density Beta with parameters $(a_{12}-a_{22})\sigma^{-2}$ and $(a_{21}-a_{11})\sigma^{-2}$.
\end{proposition}
Observe that for the deterministic system associated to~\eqref{eq:stochastic_replicator_1} the unique ergodic probability measures that support strict subsets of the community are delta distribution, namely: $\bf \delta_1$ with all the mass at the first population, and $\bf \delta_2$ with all the mass at the second population. It is direct to show that the average positivity of the invasion rates (Assumption \ref{ass:positive_Lyapunov}) holds in this case. 

Assume that we choose $0<\delta < \min \{\sigma^2, a_{12}-a_{22}, a_{21}-a_{11}\}$. For an arbitrary probability measure $\mu$ on the simplex, define $\Phi (\xx) =\tilde A \xx$, $\Psi (\xx)= \Sigma$ and
\begin{align}\label{eq:potentials}
\Upsilon (\xx, \mu)&=\Bigg ( \delta \int y_1 \mu (dy_1, dy_2)\left (\begin{matrix}
0 & 1 \\
-1 & 0
\end{matrix} \right )\Bigg )\xx,
\end{align}
and the associated non-linear replicator:
\begin{align} \label{eq:stochastic_replicator_2}
d\Xi_t =  \proj \Phi (\Xi_t) dt + \proj \Upsilon (\Xi_t, \mu_t) dt+\proj \Psi (\Xi_t)d\WW_t. 
\end{align}

\begin{theorem} \label{th:simple_theorem}
Equation \eqref{eq:stochastic_replicator_2} admits a unique invariant distribution, namely the $\text{Beta}(s, 1-s)$ distribution, where $s= \dfrac{a_{12}-a_{22}}{\sigma^2-\delta}$.
\end{theorem}
The proof of the above theorem will be transparent as soon as we prove its extension to higher dimensions (and the computation of the invariant measure will follow after some simple algebraic manipulation). 

Before starting the next Lemma, observe that if a payoff matrix $A$ satisfies \textbf{[C1]}, for every $\xx= (x_1, \ldots, x_d)^\top \in \mathcal{T}$  we have:
\begin{align*}
\xx^\top A \xx &= \dfrac{1}{2} \sum_{i,j}(a_{ij}+a_{ji})x_ix_j = \dfrac{1}{2}\sum_{i,j} (a_{ij}+a_{ji}-a_{ii}-a_{jj})x_ix_j \\
&= \dfrac{1}{2}\sum_{i\neq j} (a_{ij}+a_{ji}-a_{ii}-a_{jj})x_ix_j \\
&= \dfrac{1}{2}\sigma^2 \sum_{i=1}^d x_i (x_1+ x_2 + \ldots + x_{i-1}+x_{i+1} + \ldots +x_d)
= - \dfrac{1}{2}\sigma^2 \sum_{i=1}^d x_i^2.
\end{align*}
In particular, if $\xx \in \mathcal{T}$ is not the null vector, $A \xx$ cannot be a constant vector. We deduce that the equation $ A' \alpha= \text{ constant vector }$ has at most one solution summing up to $1$. 

\begin{lemma} \label{lemma:tolerance}
Assume that $A$ satisfies \textbf{[C1]} and \textbf{[C2]}. Then there exists $\delta > 0$ such that for any skew-symmetric matrix $E$ with $\Vert E \Vert < \delta$ the matrix $\hat A:= A+E$ satisfies \textbf{[C2]}. 
\end{lemma}

\begin{pf}
First, observe that for $A$ and $E$ as in the claim, the matrix $A+E$ satisfies \textbf{[C1]} as well. Next, observe that we can assume that $E$ has the block-form:
\begin{align*}
E= \delta
\left(
    \begin{array}{c;{2pt/2pt}c}
        E_2 & {\bf 0}_{2\times(d-2)} \\ \hdashline[2pt/2pt]
        {\bf 0}_{(d-2)\times 2} & {\bf 0} 
    \end{array}
\right),\quad\text{ with }\quad E_2=\begin{pmatrix}
0 & 1\\ -1 &0
\end{pmatrix},
\end{align*}
since, if the result is true for this kind of perturbation, applying repeatedly the lemma we arrive (via permutations of rows and columns) at the desired result.

We aim to find a strictly positive vector $\hat \alpha = (\hat \alpha_1, \ldots, \hat \alpha_d)^\top$ such that $\hat A ' \hat \alpha$ is a constant vector. Let $\alpha:= (\alpha_1, \ldots, \alpha_d)^\top \in \R{d}$ and $c \in \mathbb{R}$ be the vector and the real number ensured by~\textbf{[C2]} associated to the matrix $A$. For $i=1, 2, \ldots, d-1$, write $\hat \alpha_i:= \alpha_i + \varepsilon_i$, and $\hat \alpha_d= \alpha_d - \sum_{j=1}^{d-1} \varepsilon_j$. Then for a certain constant $\hat c$, for $i=1, \ldots, d$:
\begin{align*}
\hat c&= \sum_{j=1}^d (\hat a_{ij}- \sigma^2/2 )\hat \alpha_j\\
&= \sum_{j=1}^{d-1} ( a_{ij} + \delta e_{ij}- \sigma^2/2 ) (\alpha_j+ \varepsilon_j) + (a_{id} + \delta e_{id} - \sigma^2/2) (\alpha_d - \sum_{k=1}^{d-1}\varepsilon_k), 
\end{align*}
must hold. Then:
\begin{align*}
\hat c & = c + \delta (\alpha_2+ \varepsilon_2)+  \sum_{j=1}^{d-1} (a_{1j}-a_{1d}) \varepsilon_j\\
&= c - \delta (\alpha_1+ \varepsilon_1) + \sum_{j=1}^{d-1} (a_{2j}-a_{2d}) \varepsilon_j \\
&= c+ \sum_{j=1}^{d-1} (a_{ij}-a_{id}) \varepsilon_j, \quad i=3, 4, \ldots, d,
\end{align*}
and consequently:
\begin{align*}
\delta (\alpha_2+ \varepsilon_2)+ \sum_{j=1}^{d-1} (a_{1j}-a_{1d}) \varepsilon_j & = - \delta (\alpha_1+ \varepsilon_1) + \sum_{j=1}^{d-1} (a_{2j}-a_{2d} )\varepsilon_j \\
&= \sum_{j=1}^{d-1} (a_{ij}-a_{id}) \varepsilon_j,
\end{align*}
for $i=3, 4, \ldots, d$. This reduces to a linear system of $d-1$ equations for the $d-1$ unknowns $\varepsilon_1, \ldots, \varepsilon_{d-1}$. Let $H^{(\delta)}$ be the matrix representing the above system. Then ${\bf \varepsilon}= (\varepsilon_1, \ldots, \varepsilon_{d-1})^\top$ satisfies the system $\bf (S_{\delta})$:
\begin{align} \label{eq:linear_system}
H^{(\delta)} {\varepsilon^{(\delta)}} = \begin{pmatrix}
\delta (\alpha_1+\alpha_2) &
\delta \alpha_2 &
\delta \alpha_2 &
\cdots &
\delta \alpha_2
\end{pmatrix}^\top.
\end{align}
Since the equation $A' \alpha=\text{constant}$ is uniquely solvable for $\alpha$ strictly positive, summing up to $1$ and satisfying~\textbf{[C1]}, the equation $H^{(0)} \varepsilon'=0$ has a unique solution, and hence $H^{(0)}$ is non singular. On the other side, it is clear that $(H^{(\delta)}: \delta \in [0,1])$ is a Lipschitz continuous family of matrices (with respect to a fixed reference norm in ${\R{d}} \otimes \R{d}$, say the $L^2$ subordinate matrix norm), and thus there exists $1 \ge  \overline \delta> 0 $ such that for $\delta < \overline \delta$ there exists a unique solution of $\bf (S_{\delta})$. Let $(\varepsilon^{(\delta)}: \delta \in [0, \overline \delta])$ be the family of solutions, and let $\ell_{\delta}$ be the right hand side of $\bf (S_{\delta})$. Then for $\delta_1, \delta_2 \in [0, \overline \delta]$, and some constant $C>0$ that can change from line to line:
\begin{align*}
\Vert \varepsilon^{(\delta_1)}- \varepsilon^{(\delta_2)}\Vert & \le \Vert (H^{(\delta_1)})^{-1} (\ell_{\delta_1}- \ell_{\delta_2}) \Vert + \Vert ((H^{(\delta_1)})^{-1}- (H^{(\delta_2)})^{-1} )\ell_{\delta_2} \Vert \\
& \le C \Vert (H^{(\delta_1)})^{-1} \Vert \Vert \ell_{\delta_1}- \ell_{\delta_2} \Vert + \Vert (H^{(\delta_1)})^{-1}- (H^{(\delta_2)})^{-1} \Vert \Vert \ell_{\delta_2} \Vert \\
& \le C\Vert (H^{(\delta_1)}) ^{-1}\Vert \vert \delta_1-\delta_2\vert 
\\
&\qquad+ C \Vert (H^{(\delta_1)})^{-1} \Vert \Vert (H^{(\delta_2)})^{-1} \Vert \Vert H ^{(\delta_1)}- H^{(\delta_2)}\Vert. 
\end{align*}
As mentioned earlier, $(H^{(\delta)}: \delta \in [0, \overline\delta])$ is a Lipschitz family; in particular, it has bounded norm and, consequently, the family of inverses is bounded as well (this is so since the coefficients of $H^{(\delta)}$ are continuous in $\delta$, and consequently they are bounded on $[0, \overline\delta]$; the determinant across the family, being the determinant a continuous functions of the coefficients, is a continuous function of a parameter ranging over a compact interval, and thus is bounded; again, since $[0,\overline \delta]$ is compact,  the determinants are bounded away from zero; Cramer's rule gives the conclusion).  Again, by Lipschitz continuity, the second term is bounded by $C' \vert \delta_1- \delta_2\vert$. We deduce that $(\varepsilon^{(\delta)}: \delta \in [0, \overline \delta])$ is a Lipschitz family. So, we can find $\hat \delta < \overline \delta $ such that for $\delta < \hat \delta$ the vector $\hat \alpha$ is still strictly positive. Since, by definition, $\sum_{j=1}^d\hat \alpha_i= 1$, this proves our claim.
\end{pf}

As a last ingredient in the proof of the main Theorem of this section, we need an estimate of Wasserstein's distance between Dirichlet distributions in terms of the differences between their parameters. Hereafter, for a positive vector ${\bf a}
= (a_1, \ldots, a_d)$ with $\sum_{i=1}^d a_i=1$, the notation $D_{\bf a}$ stands for the absolutely continuous measure with $\text{Dirichlet}({\bf a})$ density. 

\begin{lemma} \label{lemma:Dirichlet_Lipschitz}
For $D_{\bf a}$ and $D_{\bf b}$ as above:
\begin{align*}
\Was_1 (D_{\bf a}, D_{\bf b}) \le \sum_{i=1}^{d-1} 2 (d-i) \vert a_i-b_i\vert.
\end{align*}
\end{lemma}

We could not find this result in the available literature, so we provide a full proof in the Appendix.

\begin{theorem}
Let $A$ be a payoff matrix satisfying \textbf{[C1]} and \textbf{[C2]}, and $F: \mathcal{P}(\simp{d}) \mapsto \R{d} \otimes \R{d}$ be a $\Was_1$-Lipschitz function with range included in the set of skew-symmetric matrices. Let $\Phi (\xx) =\tilde A \xx$, $\Psi (\xx)= \Sigma$, and $\Upsilon (\xx, \mu)=  \delta F (\mu) \xx$. Then there exists $ \delta_0 > 0$ such that for every $\delta < \delta_0$ the stochastic McKean-Vlasov replicator with the previous parameter admits a unique invariant probability measure putting no mass on the boundary of the simplex.
\label{th:perturbation}
\end{theorem}

\begin{pf}
For $A$ as above, let $\overline \delta$ be the tolerance guaranteed by Lemma \ref{lemma:tolerance}, and fix $\delta < \overline \delta$. For fixed $\mu \in \mathcal{P}(\interior(\simp{d}))$, let $(\XX^{(\mu, \delta)}_t: t \ge 0)$ be the unique solution of the  ``frozen'' replicator:
\begin{align} 
\label{eq:frozen_replicator}
d\XX^{(\mu, \delta)}_t = \proj \Phi (\XX^{(\mu, \delta)}_t) dt + \proj \Upsilon (\XX^{(\mu, \delta)}_t, \mu) dt + \proj \Psi (\XX^{(\mu, \delta)}_t)d\WW_t. 
\end{align}
By Theorem 3.6 of \cite{hofbauerimhof2009}, there exists a vector $\alpha= \alpha (\mu)$ of strictly positive entries summing up to $1$ such that $D_{\alpha (\mu)}$ is invariant for $\XX^{(\mu, \delta)}$. Moreover, by Corollary 3.8 of \cite{hofbauerimhof2009}, this is the unique invariant distribution that does not charge the boundary. Let $T_\mu= D_{\alpha(\mu)}$. Moreover, since $F$ is Lipschitz, from the proof of Lemma \ref{lemma:tolerance}, it transpires that there exists a constant $C> 0$ such that: 
\begin{align*}
\Vert \alpha_{\mu}- \alpha_{\nu}\Vert_1 \le \delta C \Was_1 (\mu, \nu),
\end{align*}
for every $\mu, \nu \in \mathcal{P} (\simp{d})$. But then, by Lemma \ref{lemma:Dirichlet_Lipschitz}:
\begin{align*}
\Was_1 (T_\mu, T_\nu) \le 2(d-1) \delta C \Was_1 (\mu, \nu).
\end{align*}
We deduce that for $\delta < \min \{\overline\delta, (2(d-1)C)^{-1}\}: = \delta_0$, the map $\mu \mapsto T_\mu$ is a contraction. By the Contraction Mapping Theorem (see e.g. Theorem 5.7, page 138 of \cite{brezis2010}), there exists a unique $\overline \mu\in \mathcal{P} (\simp{d})$ such that $\overline \mu= T_{\overline \mu}$. Clearly, $\overline \mu$ is invariant for the McKean-Vlasov replicator. Since the range of $T$ is a subset of the Dirichlet family, the last claim follows.   
\end{pf}

\section{Simulations} \label{sec:simulations}

In Section \ref{section:invariant}, for the case of the McKean-Vlasov replicator with linear interaction function and (small) Lipschitz dependence on the instantaneous measure (see Theorem~\ref{th:perturbation}), we have stated the existence of a unique invariant probability measure putting no mass on the boundary. In this section, via numerical simulations, we investigate whether this invariant probability attracts the solutions of the system at large times. To that aim, we use the well-known Euler-Maruyama scheme to simulate the associated $N$-replicator system, with step size $h=0.01$. 

Our objective is to study the behaviour of the solutions to~\eqref{eq:stochastic_replicator_2} as $t\rightarrow\infty$, i.e. the long term solutions to
\begin{equation*}
dX^{(N;i)}_t =  \proj \Phi (X^{(N; i)}_t) dt + \dfrac{1}{N}\sum_{j=1}^N\proj \Upsilon (X^{(N; i)}_t, X^{(N;j)}_t) dt + \proj \Psi (X^{(N; i)}_t) dW_t,
\end{equation*}
where $X_t^{(N;i)} = (x_1^{(N;i)}(t),x_2^{(N;i)}(t))$ for $N$ large. Here $\Phi,\Psi$ and $\Upsilon$, correspond to those defined in \eqref{eq:potentials}. In particular, we simulated the following sets of parameters
\begin{equation}
\tag{PS1}
A = \begin{pmatrix}0.5 & 1\\1 & 0.5\end{pmatrix},\qquad \sigma=1,\qquad \delta = 0.05,\qquad s=0.5263,
\end{equation}
and
\begin{equation}
\tag{PS2}
A = \begin{pmatrix}0.6 & 0.9\\1 & 0.4\end{pmatrix},\qquad \sigma=0.9487,\qquad \delta = 0.04,\qquad s=0.5814,
\end{equation}
For each of the two sets of parameters, to test the effect of the initial distribution, we perform the simulations for initial conditions uniformly scattered in the unit interval, and also for initial conditions uniformly scattered in the subset $[0.2,0.4]$. We called those distributions \textit{uniform} and \textit{lump} initial conditions (UIC and LIC for short). For each of the four experiments described, we approximate the probability density function from a normalized histogram of the first coordinate of each of the $N$-replicator $\{x_1^{(N;i)}\}_{i=1}^N$. This gives us some qualitative understanding of the time convergence of the probability density. 

\begin{figure}
\centering
\begin{minipage}{0.9\textwidth}
\includegraphics[width=\textwidth]{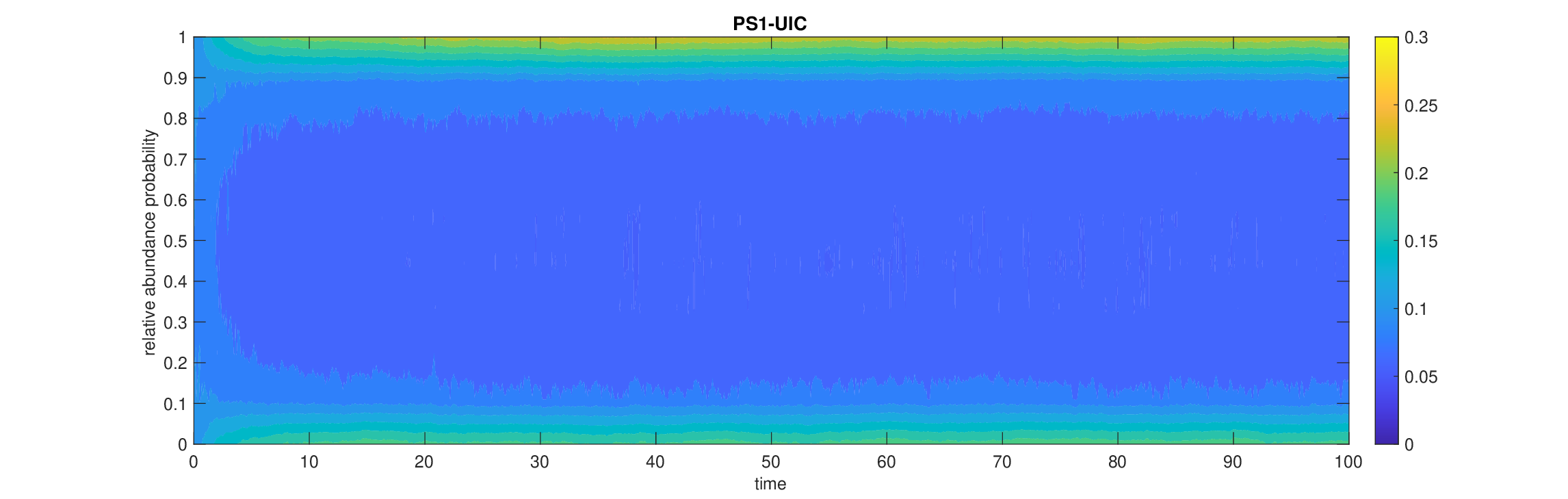}
\end{minipage}
\hfill
\begin{minipage}{0.5\textwidth}
\includegraphics[width=\textwidth]{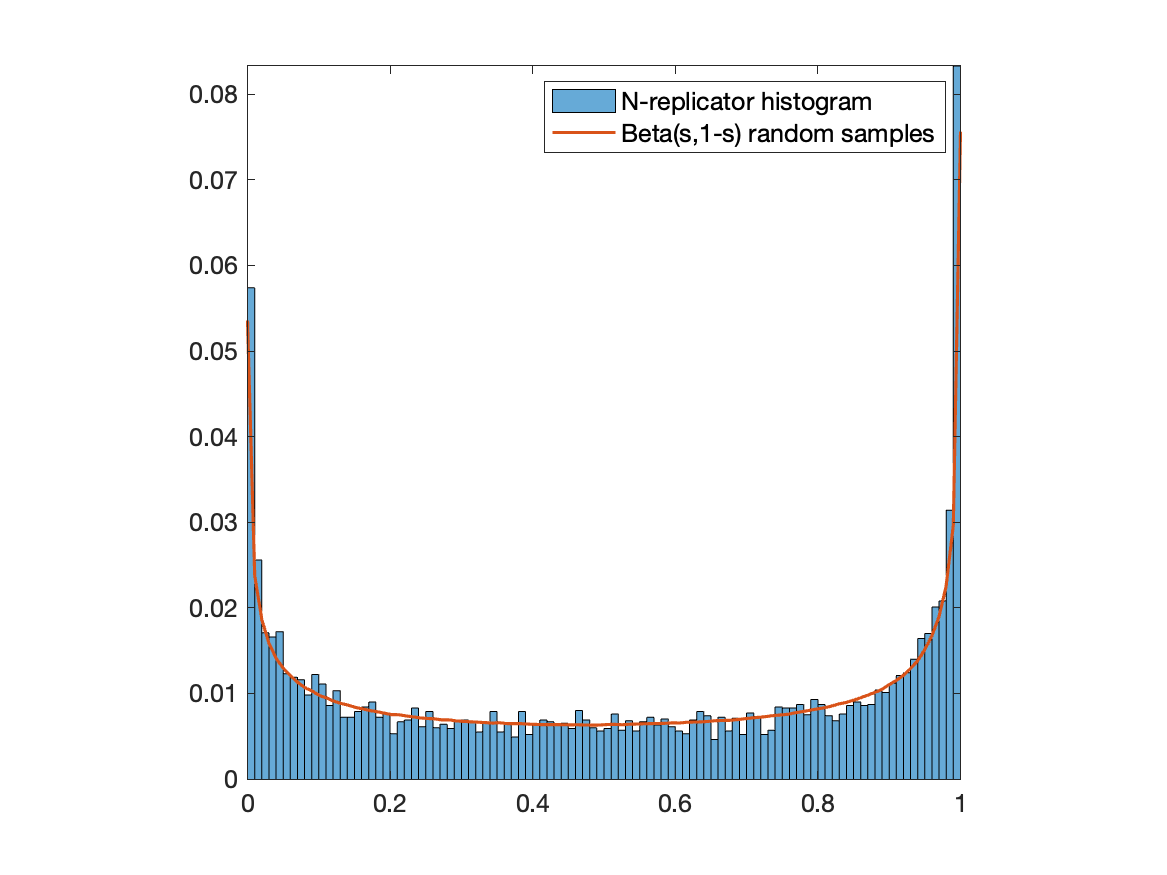}
\end{minipage}

\caption{Simulations for the $(\simp{2})^{10000}$ interacting particle system. For the $N$-replicator system, at each time step, we calculate the empirical histogram from the first coordinate of each replicator under probability normalization. The top panel shows the time evolution of those empirical relative abundance probabilities for a set of parameters and initial conditions PS1-UIC. The bottom panel shows the histogram for $N$-replicator at $t=100$ (bins) and for a sample of $10^7$ random numbers drawn for the theoretical Beta(s,1-s) (red line).}
   \label{fig:2d_simplex}
\end{figure}

The contour plot shown in Figure~\ref{fig:2d_simplex} (top panel) corresponds to the time evolution of such histograms for PS1 and UIC. Results for LIC and PS2 can be found in the appendix Figure~\ref{fig:a2d_simplex}. We notice that the histograms stabilize after approximately 3000-time steps ($t=30$), with small fluctuations due to finite-size effects. At the bottom panel of Figure~\ref{fig:2d_simplex} we see that the final empirical stationary pattern is similar to the Beta distribution stated at that Theorem~\ref{th:simple_theorem}, and it seems to not depend on the initial distributions UIC or LIC (see appendix Figure~\ref{fig:a2d_simplex}). If such convergence holds true, then as the time $t$ goes by, the position of the first coordinate of each of the $N$ interacting replicators would constitute a sample from the unique Beta distribution. In that case, the empirical mean $\frac{1}{N}\sum_{i=1}^N x_1^{(N;i)}$ would be converging towards the mean of the limit distribution which is, in this case, $s/(s+(1-s))=s$. To test this conjecture we calculate - for each parameter set and initial distribution - at each time step the empirical mean and calculate the relative percentage error with respect to $s$. For any of the numerical experiments performed, the results showed a relatively fast convergence of the empirical mean towards $s$, with small fluctuations around $s$ smaller than 2\%.



\begin{figure}[!ht]
 \centering
\includegraphics[width=0.44\textwidth]{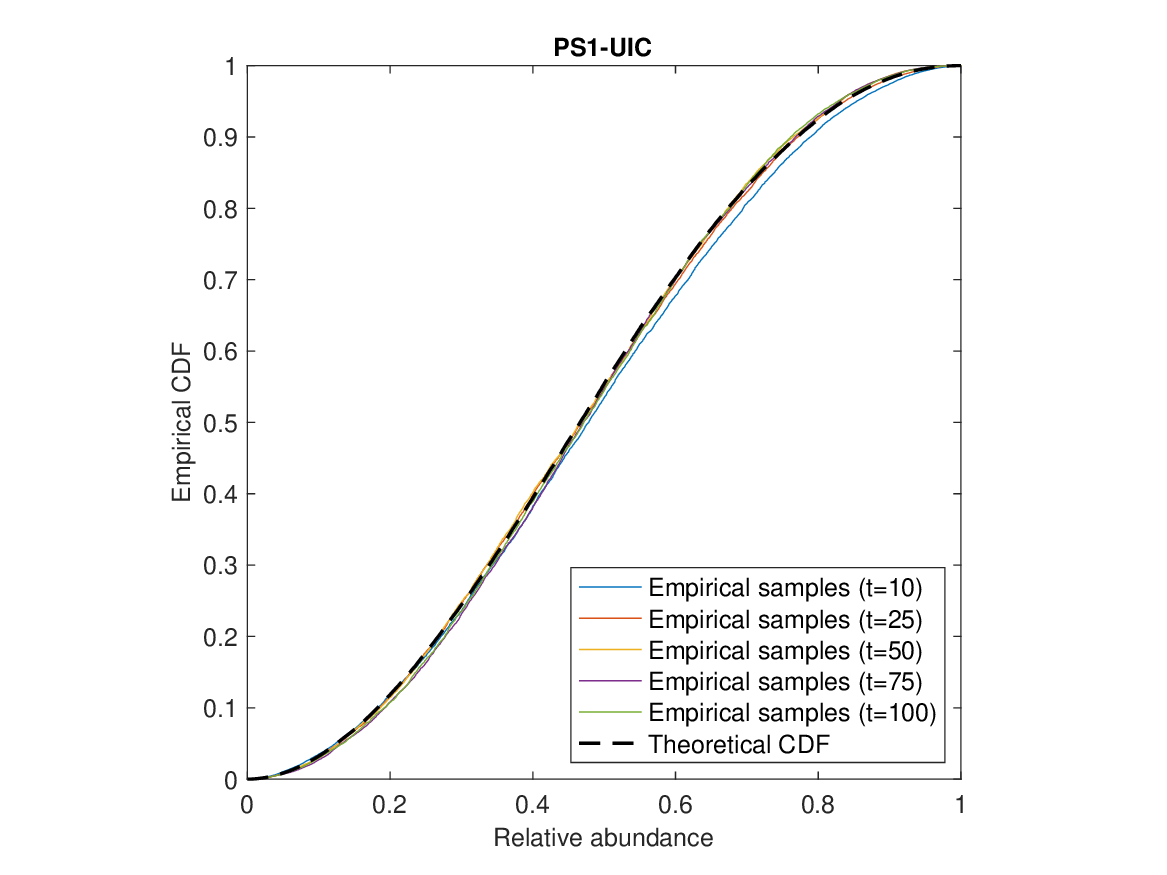}
\hfill
\includegraphics[width=0.44\textwidth]{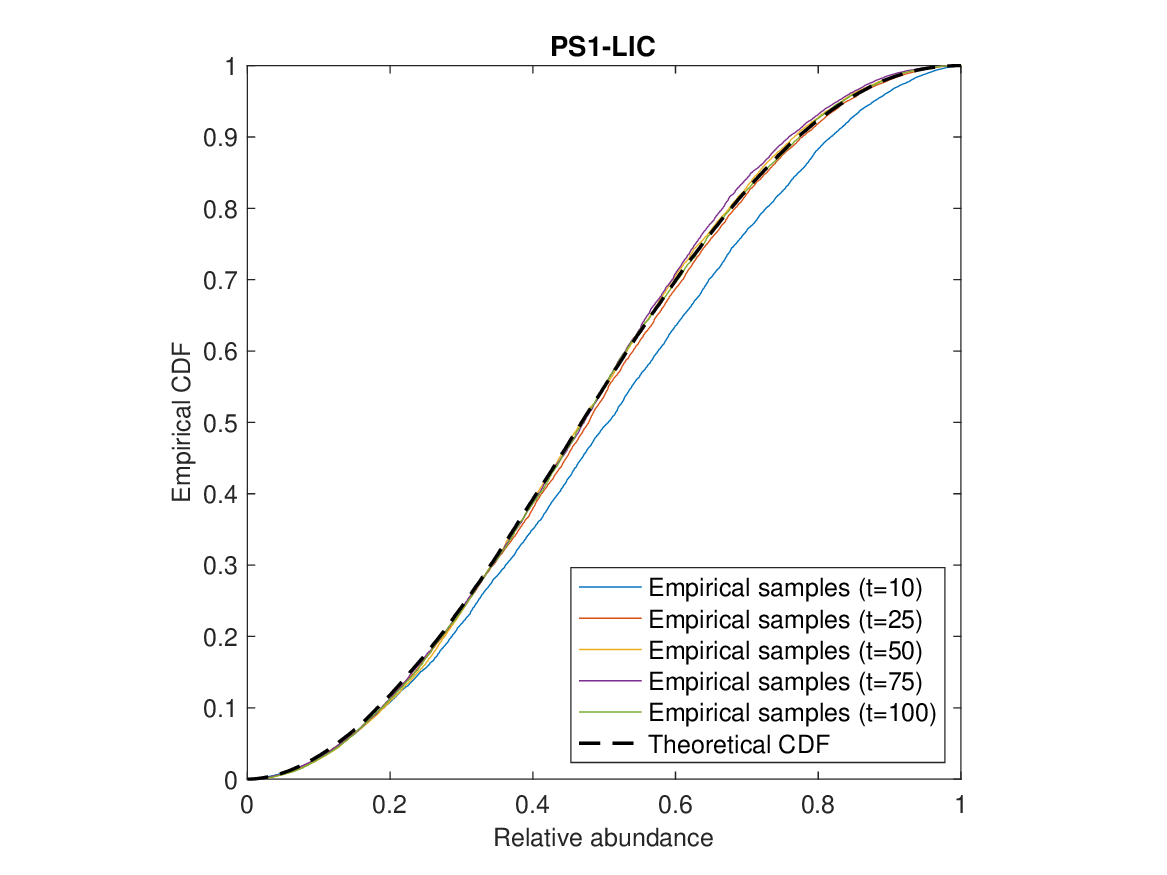}

\includegraphics[width=0.44\textwidth]{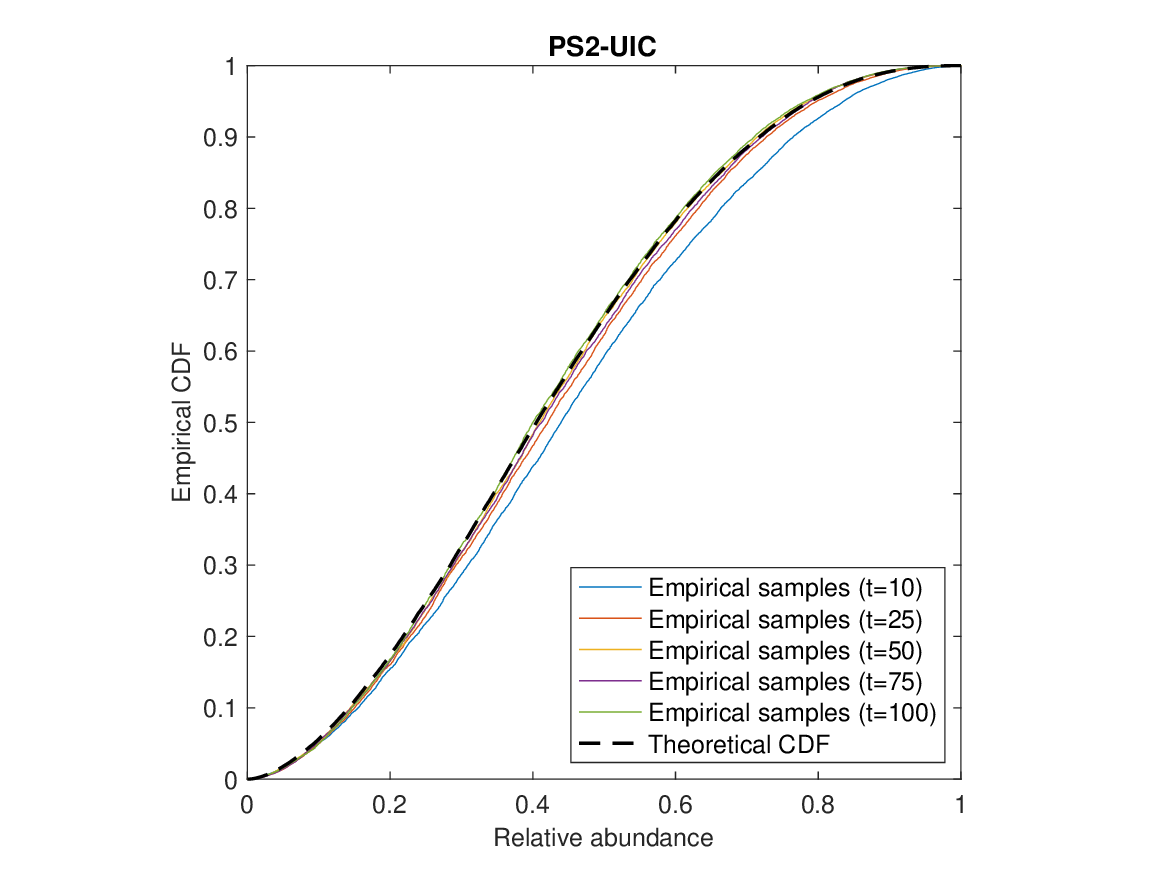}
\hfill
\includegraphics[width=0.44\textwidth]{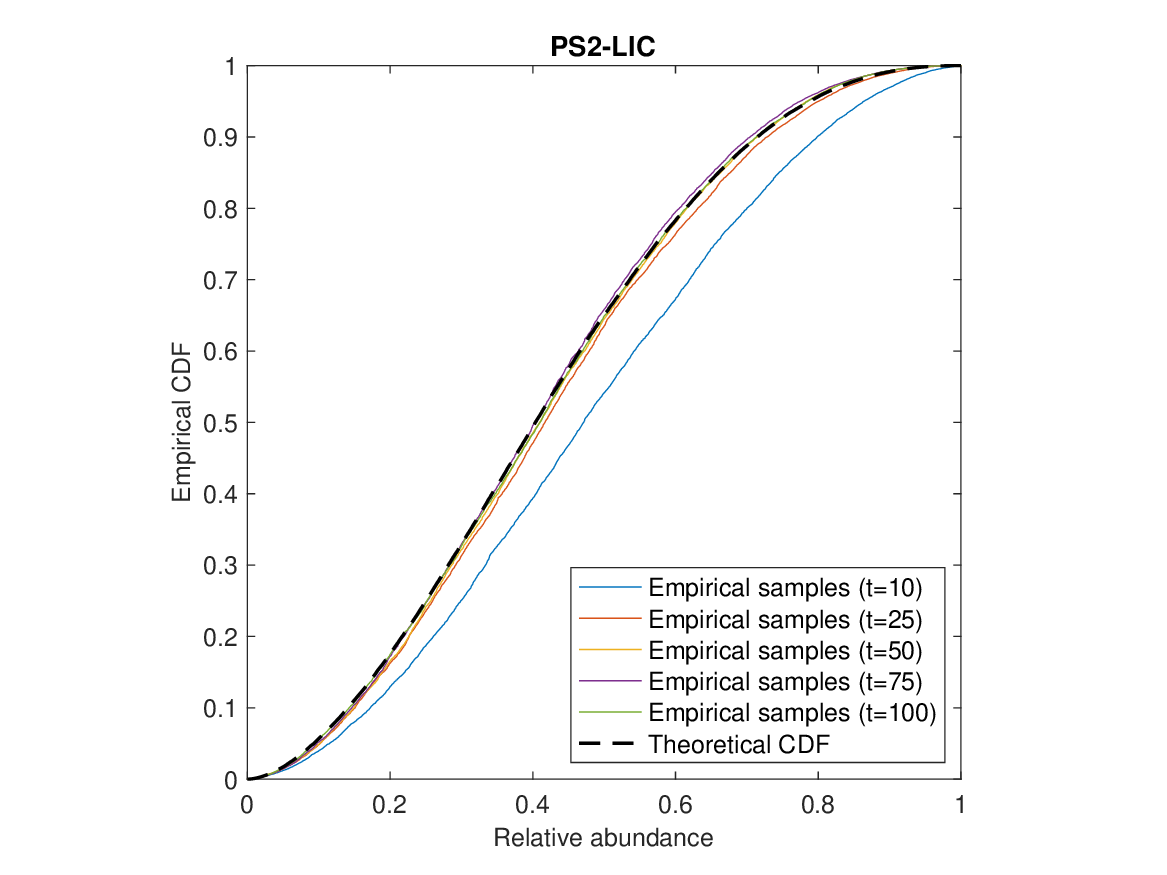}
    
\caption{Empirical cumulative distribution function of the first coordinate of the $N$-replicator system. At the top (left) panels, we show the first parameter set (uniform initial distribution), and in the bottom (right) panels we show the second parameter set (lump initial distribution). In all plots, we show the ECDF for times $t=10,25,50,75$ and $t=100$, and the dashed black line is the CDF associated with the theoretical Beta(s,1-s) distribution. Again, for times larger than $25$, the shape of the empirical and theoretical distribution are qualitatively similar supporting the conjecture of ergodicity.}
   \label{fig:2d_cdf}
\end{figure}

To test that the convergence is not only restricted to the empirical mean but also to the empirical distribution given by the simulated $N$-replicator system, we calculate the empirical cumulative distribution function (ECDF) for each set of parameters and initial distributions. In Figure~\ref{fig:2d_cdf} we plot the ECDF for different times and each numerical experiment. In all cases, the shape of the empirical CDF and the theoretical CDF given by Theorem~\ref{th:simple_theorem}, are relatively similar. To further study this hypothesis, we performed at each time step $t_n$ larger than $t=30$ an Anderson-Darling test to see whether the sample vector $x_1^{(N;i)}(t_n)$ is drawn from the theoretical Beta(s,1-s) distribution. For a 1\% of significance level, we see that there is no statistical evidence to reject the null hypothesis for 64.91\%, 72.47\%, 76.32\%, and 80.31\% of the times for the four PS1-UIC, PS1-LIC, PS2-UIC and PS2-LIC respectively.

As a consequence, it seems reasonable to conjecture that ergodicity must hold true. We perform an additional statistical analysis, but now consider true independent trajectories of $N$-replicators at a significantly lower $N$ than the previously used ones. To that end, we now simulate $10^3$ times the $500$-replicator dynamics with UIC for both sets parameters, saving the trajectory of only one randomly selected replicator each time. We report a goodness-of-fit test regarding the theoretical stationary Beta distribution specified in Theorem \ref{th:simple_theorem} at times $t=40$, and $t=50$. In all cases, we will see now that the laws of the trajectories are quite close to the corresponding theoretical stationary Beta distribution (for more details see also appendix Figure \ref{fig:NewBeta}).

The implemented test is based on an $L^2$-type test statistic proposed in \cite[Theorem 1 - Corollary 1]{ebner2021new}: based on a random sample $Y_1, ..., Y_n$, their test statistic is given by:
\[
T_n(\hat{a}_n,\hat{b}_n)=n\int_{0}^{1}\left|\frac{1}{n}\sum_{j=1}^{n}\left((\hat{a}_n+\hat{b}_n)Y_j-\hat{a}_n\right)\mathbf{1}_{\{Y_j\geq t\}}-\frac{t^{\hat{a}_n}(1-t)^{\hat{b}_n}}{B(\hat{a}_n,\hat{b}_n)}\right|^{2}dt,
\]
where $B(\cdot,\cdot)$ is the Beta function, and $\hat{a}_n$ and $\hat{b}_n$ are consistent estimators of the Beta parameters, such as moment estimators or maximum likelihood estimators (MLEs). However, since in our simulations, the stationary distribution is completely specified by the resulting parameterization given in Theorem \ref{th:simple_theorem}, we can introduce such values, $a$ and $b$, instead of its estimates. Thus, for our null hypothesis,  $H_0$: ``the law of the  first coordinate trajectories of the one-dimensional McKean-Vlasov replicator \eqref{eq:stochastic_replicator_2} follows a $Beta(a,b)$ for time $t$'',  our test statistic will become:
\begin{equation}\label{stat-test}
T_n=T_n(a,b)=n\int_{0}^{1}\left|\frac{1}{n}\sum_{j=1}^{n}\left((a+b)Y_j-a\right)\mathbf{1}_{\{Y_j\geq t\}}-\frac{t^{a}(1-t)^{b}}{B(a,b)}\right|^{2}dt,
\end{equation}
where in our case $a=s$ and $b=1-a$. The critical region to reject the $H_0$ is delimited by $1-\alpha$ quantiles of the distribution of the above test statistic under $H_0$. Since $T_n$ is a distance, the larger the observed values of $T_n$ are, the further we will be from the null hypothesis. To obtain an approximation of the distribution of $T_n$ under $H_0$, we construct the density of  $B$ values of $T_n$, each one obtained by an independent simulated random sample from a $Beta(a,b)$. In our case, we used $B=5000$. We found that the $90\%$ quantiles are $0.1056$ (PS1) and $0.1119$ (PS2). For $t=40$ the computed tests statistics are $T_n=0.0175$ (PS1) and $T_n=0.0046$ (PS2). For $t=50$ the computed tests statistics are $T_n=0.0817$ (PS1) and $T_n=0.0492$ (PS2). In all cases, the null hypothesis is not rejected. The corresponding \textbf{R}-code can be found at \url{https://github.com/leonardo-videla/beta_test.git}.


Given these results, it seems reasonable to conjecture that ergodicity must hold, at least in this case. We have no proof yet, and thus this remains an open question.

\section{Discussion and outlook} \label{sec:conclusions}

This paper advances the modeling of biological phenomena, which we have proposed in previous papers using the framework of Open System Dynamics (\cite{marquet2017proportional}, \cite{Rebolledo.2019}, \cite{freilich2020reconstructing}, \cite{marquet2020species}, \cite{tejo2021coexistence}). 
We regard our contribution as a step forward in understanding how the living develops in continuous interaction with the environment, including this time the interaction within communities of living entities in what we envisioned as an ecology of primordial or early life represented by a community of replicators. In particular, we were interested in assessing if neutrality and associated fitness equivalence could emerge in this system and under which conditions. To this end, we generalized the stochastic replicator equation by adding a new component that affects the instantaneous fitness, the idea being that this component accounts for the interaction of one particular replicator population with others, alike ones, within a community in a mean-field regime.

We show that under suitable hypotheses, the asymptotic behavior of interacting replicator dynamics may be approached by independent non-markovian, McKean-Vlasov particles (Section \ref{sec:propchaos}). For the mean-field replicator aggregates, we have shown that the conditions used in \cite{benaim2008robust}, with some simple modifications, still apply to guarantee the coexistence (in the sense of persistence) of large (but finite) systems of interacting communities (Section \ref{sec:persistence}). In simple cases (Section \ref{section:invariant}), we show the existence of a unique non-trivial invariant probability measure for the McKean-Vlasov replicator. Indeed, in simple setups, we have been able to prove that the Dirichlet family is stable (in the obvious sense that a system with Dirichlet invariant still has a Dirichlet invariant when we introduce the dependence on the measure). This result, which partially extends the work \cite{hofbauerimhof2009} and \cite{imhof2005}, admits further extensions, and it is a part of our current research efforts. Likewise, regarding the ergodicity of the McKean-Vlasov replicator, we still have no concluding results (although we have shown in Section \ref{sec:simulations} that numerical experiments give strong empirical support in this direction). This is an open question we continue investigating. 

We point out some salient features of this article. The first one refers to the fact that, to the best of our knowledge, this is the first attempt to import the McKean-Vlasov machinery into the evolutionary-game-theoretic realm. Also, the study of persistence for interacting, finite aggregates of complex communities seems to be a novel feature. The use of the fixed-point argument to prove existence/uniqueness of invariant measures for McKean-Vlasov systems (see Theorem \ref{th:invariant}) also seems to be a technique that has not been previously informed in this literature. 

From the point of view of theoretical ecology, we note that the condition \textbf{[C1]} can be regarded as the expression of neutrality and fitness equivalence among types, as assumed in the neutral ecological theory (\cite{hubbell2001}); in turn, in previous work (\cite{marquet2017proportional}) we have shown that this condition is associated with Beta invariant probability measure. In light of the results of Section \ref{section:invariant}, this observation provides theoretical support for the idea that neutrality, or fitness equivalence, is associated with the emergence of ecologies from simple, replicating molecules to complex ecosystems, and is an emergent condition that is related to the persistence of the ecological system. Thus, unlike traditional models in ecology, we did not assume neutrality but neutrality naturally emerged as a condition for persistence.

There are many open directions to continue this work. The first one was alluded to above, namely: a proof for ergodicity of the McKean-Vlasov replicator is still missing. Secondly, we were able to prove persistence for the finite-size particle system under small mean-field and noise terms. The problem of proving the analogous property for the McKean-Vlasov replicator was not faced in this article, but we notice in advance that the techniques used in Section \ref{sec:persistence} (a modification of the well-known Foster-Lyapunov techniques as treated for example in the classic articles of Meyn and Tweedie \cite{meyntweedie92,meyntweedie93a,meyntweedie93b}) are unlikely to work in the non-linear case since they are strongly grounded on the Markovian nature of the objects they are meant to be applied to. Finally, another aspect that must be investigated is the extension of the model proposed here to the setup of many types of interacting replicators. More specifically, we can think of a myriad of entities of different types, say $N_i\ge 1$ replicators of type $i$ for $i=1, \ldots, M$, where $M$ is the fixed number of types. A preliminary assessment of this extended setup shows that the propagation-of-chaos (under suitable conditions) must still hold. The problem of persistence in this multi-type interacting replicator seems to be a challenging issue, and we think it deserves further investigation. 
\smallskip

\begin{small}
\noindent\textbf{Acknowledgements}
This work was funded by project FONDECYT 1200925 \emph{The emergence of ecologies through metabolic cooperation and recursive organization}, Centro de Modelamiento Matem\'atico (CMM),
Grant FB210005, BASAL funds for centers of excellence from ANID-Chile, Exploration-ANID 13220168, \emph{Biological and quantum Open System Dynamics: evolution, innovation and mathematical foundations}, FONDECYT Iniciaci\'on project number 11240158-2024 \emph{Adaptive behavior in stochastic population dynamics and non-linear Markov processes in ecoevolutionary modeling}, and FONDECYT Iniciaci\'on 11200436, \emph{Excitation and inhibition balance as a dynamical process}, and \textit{``Programa de Inserción Acad\'emica 2024 Vicerrector\'ia Acad\'emica y Prorrectoría de la Pontificia Universidad Cat\'olica de Chile''}. We thank Professor Rodrigo Plaza for his valuable help in the implementation of statistical tests. We also thank the editor and an anonymous referee who provided us with valuable comments that greatly improved our manuscript. 
\end{small}
\smallskip

\noindent\textbf{Declarations and conflict of interest}

The authors declare that have no financial or personal relationship with other people or organizations that could inappropriately influence or bias the content of this paper. Also, the authors declare that they have no conflicts of interest. 

\bibliography{./bibliography_simplex.bib}
\bibliographystyle{abbrv}

\appendix 

\section{Appendix}

\subsection{Derivation of the replicator equation from a model of population dynamics} 

Define $\bG: \mathbb{R}^d\setminus\{0\} \mapsto \simp{d}$ via:
\begin{align*}
\bG (\yy) = \Vert \yy \Vert_1^{-1} \yy,
\end{align*}
and consider the SDE:
\begin{align*}
d\YY_t = \YY_t \circ (\tilde \Phi (\bG (\YY_t))dt+ \Psi (\bG (\YY_t))d\WW_t).
\end{align*}
Here, $\YY$ is to be understood as a $\mathbb{R}^d_{+}$-valued process of abundances. The characteristic $\tilde \Phi$ and $\Psi$ are suitable vector- and matrix-valued functions, respectively. Observe that in this formulation, the instantaneous fitness rates of the populations depend on $\YY_t$ only through the relative proportions of the populations. Consider the case where $\Psi$ is a diagonal matrix (the general case is analogous). Write $\XX_t= (X^{(i)}_t)_{i=1, \ldots, d}= \bG (\YY_t) = (g_i (\YY_t))_{i=1, \ldots, d}$, and let $\tilde \phi_i, \psi_i$ the components of $\tilde \Phi$ and the diagonal elements of $\Psi$, respectively. We easily compute:
\begin{align*}
    \partial_k g_i (\yy) &= \delta_{ik} \dfrac{1}{\Vert \yy \Vert_1 } - \dfrac{y_i}{\Vert \yy \Vert^2}. \\
    \partial_{kk} g_i(\yy) &= 2 \left ( - \delta_{ik}\dfrac{1}{\Vert \yy\Vert^2} + \dfrac{y_i}{\Vert \yy \Vert^3}\right ).
\end{align*}
Thus, It\^o's formula gives:
\begin{align*}
dX^{(i)}_t &= \Bigg(\sum_{k} Y^{(k)}_t \tilde \phi_k (\XX_t) \left(\dfrac{\delta_{ik}}{\Vert \YY_t\Vert} - \dfrac{Y^{(i)}_t}{\Vert \YY_t\Vert^2}\right)  \\
& \quad + \sum_{k} (\psi_k (\XX_t) Y^{(k)}_t)^2 \left ( \dfrac{Y^{(i)}_t}{\Vert \YY_t\Vert^3} - \dfrac{\delta_{ik}}{\Vert \YY_t\Vert^2}\right) \Bigg) dt \\
& \quad + \sum_{k} Y^{(k)}_t \tilde \psi_k (\XX_t) \left(\dfrac{\delta_{ik}}{\Vert \YY_t\Vert} - \dfrac{Y^{(i)}_t}{\Vert \YY_t\Vert^2}\right)dW^{(i)}_t\\
&=\Bigg(\sum_{k} X^{(k)}_t \tilde \phi_k (\XX_t) \left(\delta_{ik} - X^{(i)}_t \right)  \\
& \quad + \sum_{k} (\psi_k (\XX_t) X^{(k)}_t)^2 \left ( X^{(i)}_t - \delta_{ik}\right)\Bigg) dt \\
& \quad + \sum_{k} X^{(k)}_t \tilde \psi_k (\XX_t) \left(\delta_{ik} - X^{(i)}_t \right)dW^{(i)}_t.
\end{align*}
Hence, we see that $\XX$ satisfies the equation
\eqref{eq:eq00} with: 
\begin{align*}
\Phi (\xx) = \tilde \Phi (\xx)-\Psi\Psi^\top(\xx)\xx, 
\end{align*}

\subsection{Proof of Proposition~\ref{prop:pathwise_uniqueness}}

\noindent\textbf{Proposition 2:} Under Assumptions \ref{ass:non_degenerate_noise} and \ref{ass:regularity_coeff}, equation \eqref{eq:mckean_vlasov} has a unique pathwise (hence strong) solution.   
\bigskip

\label{app:uniqueness}
    \begin{pf}
    Since ${\bf b}=\proj\Phi+ \proj \Upsilon$ is Lipschitz, in particular, it is Dini-continuous. On the other hand, $\simp{d}$ is compact, and consequently the drift and diffusion coefficients are bounded. Finally, Assumption 4.2 implies that for every $\xx\in \simp{d}$, we have $\Vert \proj \Psi(\xx)^\top  \yy \Vert \neq 0$ whenever $\yy \in \mathcal{T}, \yy \neq 0$. But then:
    \begin{align*}
        \inf_{\xx \in \simp{d}, \yy \in \mathcal{T}, \Vert \yy \Vert =1 } \yy^ \top  \proj \Psi(\xx)\proj \Psi(\xx)^\top \yy > 0,
    \end{align*}
    since the expression inside the infimum is a continuous function of $\xx, \yy$ and it is taken on a compact set. Theorem 1 of \cite{veretennikov2021pathwise} yields the claim. 
\end{pf}

\subsection{Proof of Lemma~\ref{lemma:non_vanishing}}
\label{app:lemma1}

\noindent\textbf{Lemma 1:} For every $\Xi_0 \in \interior (\simp {d})$:
\begin{align*}
    \mathbb{P}^{(MV)}_{\Xi_0}\left( \min_{i=1, \ldots, d} \Xi_t^{(i)} > 0  \text{ for all } t \ge 0 \right )=1. 
\end{align*}
An analogous result holds for the $N$- particle system, namely: for every $\vec \xx_0 \in \interior ((\simp{d})^N)$,
\begin{align*}
    \mathbb{P}^{(N)}_{\vec \xx_0}\left( \min_{j=1, \ldots, N} \Vert \XX^{(N; j)}_t \Vert > 0  \text{ for all } t \ge 0 \right )=1. 
\end{align*}
In words: with probability $1$, in finite time there is no absorption on the boundary of $\simp{d}$ (for the McKean-Vlasov replicator) nor on the boundary of $(\simp{d})^{N}$ (for the $N$-particle system).
\bigskip

\begin{pf}
We prove the result for the McKean-Vlasov system, since the proof for the $N$-replicator is completely analogous. For any $k \in \mathbb{N} $, set:
\begin{align*}
V_k:= \{\xx \in \simp {d}: \min_{i} x_i < e^{-k}\},
\end{align*}
and consider the function $U_k: \simp{d} \setminus V_k \mapsto \mathbb{R}^{+}$ given by:
\begin{align*}
    U_k (\xx):= -\sum_{i} \ln (x_i).
\end{align*}
Then $U_k$ is uniformly continuous for each $k$, and indeed it is a $\mathcal{C}^{\infty}$ function in its domain.  For $\mu \in \mathcal{P}(\simp{d})$ consider the operator $L_{\mu}$ acting on a smooth function $\varphi \in \mathcal{C}^2 (\simp{d})$ as:
\begin{align}\label{eq:operator_non_linear_mckean_vlasov}
	L_\mu \varphi (\xx)=  \inner{{\bf b} (\xx, \mu)} {\nabla \varphi (\xx)} + \dfrac{1}{2}{\bf Trace} (\proj \Psi (\xx)^\top \nabla \nabla^\top \varphi (\xx) \proj \Psi (\xx) ),
\end{align} 
where $\mu_t$ solves the non-linear Fokker-Planck equation:
\begin{align}\label{eq:PDE_1}
	\dfrac{d}{dt} \inner{\mu_t}{\varphi} = \inner{\mu_t}{L_{\mu_t}\varphi},
\end{align}
for every $\varphi \in \mathcal{C}^2 (\simp{d})$ (\cite{meleard96}, \cite{chaintron2022}, \cite{sznitman91}). 

So, for $\xx \in \simp{d} \setminus V_k$ and $\mu \in \mathcal{P}(\simp{d})$, we can compute:
\begin{equation*}
    L_\mu U_k (\xx) = -\sum_i \big(\Phi_i(\xx) - \langle \xx, \Phi (\xx)\rangle\big) - \sum_{i}  \big (\Upsilon_i(\xx, \mu) - \langle \xx,  \Upsilon (\xx, \mu)\rangle \big)  + \dfrac{1}{2} \sum_i \Trace(\A_\Psi \A_\Psi^\top) (\xx), 
\end{equation*}
and we observe that $L_\mu U_k$ is, for each $k$, a uniformly continuous function on $\simp{d}$. Moreover, there exists an absolute constant, say $M$, such that $L_\mu U_k < M$ for every $k$. Let:
\begin{align*}
    \tau_k:= \inf \{t \ge 0: \Xi_t \in V_k\}.
\end{align*}
Now, fix $\Xi_0 \in \interior (\simp{d})$. There exists $k_0$ such that $\Xi_0 \in \simp{d} \setminus V_k$ for every $k \ge k_0$. It is known that for $k \ge k_0$, the process $(U_k(\Xi_{t \wedge \tau_k})-U_k (\Xi_0)- \int_{0}^t L_{\Law (\Xi_{s \wedge \tau_k})}(\Xi_{s \wedge \tau_k})ds: t \ge 0)$ is a martingale. Thus, for every $N \in \mathbb{N}$:
\begin{align}\label{Uk}
    \mathbb{E}_{\Xi_0} (U_k (\Xi_{\tau_k \wedge N})) & = U_k(\Xi_0) + \mathbb{E}_{\Xi_0} \left ( \int_{0}^{\tau_k \wedge N} L_{\Law (\Xi_{s \wedge \tau_k})} U_k (\Xi_{s \wedge \tau_k}) \mathrm{d}s \right ) \nonumber \\
    & \le U_k(\Xi_0) + MN.
\end{align}
Let:
\begin{align*}
    \Omega_N:= \bigcap_{k \ge k_0} \{ \omega \in \Omega: \tau_k < N \},
\end{align*}
and define $\epsilon_N := \mathbb{P}_{\Xi_0} (\Omega_N)$. We have:
\begin{align*}
    \mathbb{E}_{\Xi_0} (U_k (\Xi_{\tau_k \wedge N})) & \ge  \mathbb{E}_{\Xi_0} (U_k(\Xi_{\tau_k}) \cdot \mathbf{1}_{\Omega_N} )\\
    &\ge k \epsilon_N,
\end{align*}
by path-continuity. Using \eqref{Uk}, we deduce that for fixed $N$ and every $k \ge k_0$:
\begin{align*}
    k \epsilon_N \le U_k(\Xi_0) + MN \le C + MN.
\end{align*}
for some constant $C$ that depends on $\Xi_0$ but not on $k$. Letting $k$ go to infinity, we obtain $\epsilon_N=0$ for every $N$, and consequently:
\begin{align*}
    \mathbb{P}_{\Xi_0} (\lim_{k} \tau_k < \infty) = \mathbb{P}_{\Xi_0} \left ( \bigcup_{N \in \mathbb{N}}\Omega_N\right ) = 0,
\end{align*}
and this proves our claim.
\end{pf}

\subsection{Proof of Lemma \ref{lemma:Dirichlet_Lipschitz}}

\noindent\textbf{Lemma 5:} Let $D_{\bf a}$ and $D_{\bf b}$ be the Dirichlet distributions with paramters ${\bf a}$ and ${\bf b}$, both summing up to $1$. Then:
\begin{align*}
\Was_1 (D_{\bf a}, D_{\bf b}) \le \sum_{i=1}^{d-1} 2 (d-i) \vert a_i-b_i\vert.
\end{align*}
\bigskip

\begin{proof}
By dominated convergence, it is direct that the map ${\bf a} \mapsto D_{\bf a}$ is weakly continuous, i.e. continuous when we endow the set of probability measures on the simplex with the topology of weak convergence of probability measures. Thus, it suffices to consider the case where ${\bf a}, {\bf b}$ have rational entries. So, for a certain natural number $N$ there exists natural numbers $m_1, \ldots, m_d$ and $n_1, \ldots, n_d$ such that for $a_i=m_i/N, b_i=n_i/N$. Put $M_i= \sum_{j \le i} m_j$, and analogously $N_i= \sum_{j \le i}n_j$. Let $(G_j: j=1, \ldots, N)$ be an i.i.d. family of $\text{Gamma} (1/N, 1)$ random variables. Put:
\begin{align*}
V_i:= \sum_{j=M_{i-1}+1}^{M_i} G_j; \quad W_i:= \sum_{j=N_{i-1}+1}^{N_i}G_j\; \quad Z := \sum_{j=1}^n G_j.
\end{align*}
Then $X=(X_i)_{i=1}^d := (V_i/Z)_{i=1}^d$ and $Y=(Y_i)_{i=1}^d := (W_i/Z)_{i=1}^d$ have law $D_{\bf a}$ and $D_{\bf b}$ respectively. We say that the $i$ variates overlap if $[M_{i-1}+1, M_i] \cap [N_{i-1}+1, N_i] \neq \emptyset$. Observe that if the $i$ variates overlap, then $\vert V_i-W_i\vert $ is bounded by the sum of two independent Gamma-distributed random variables of parameters $(\frac{\vert N_{i-1}-M_{i-1}\vert}{N}, 1)  $ and $(\frac{\vert N_{i}-M_{i}\vert}{N}, 1)$. If the $i$ variates do not overlap, then $\vert V_i-W_i\vert $ is bounded by the sum of two (not independent) Gamma-distributed random variables of parameters $(\frac{\vert N_{i-1}-M_{i-1}\vert}{N}, 1)  $ and $(\frac{\vert N_{i}-M_{i}\vert}{N}, 1)$. Denote by $R_i$ and $S_i$ these random variables. We have:
\begin{align*}
\E (\vert X_i-Y_i\vert) &\le \E (R_i/Z) + \E(S_i/Z)\\
&= \dfrac{\vert N_{i-1}-M_{i-1}\vert}{N} + \dfrac{\vert N_{i}-M_{i}\vert}{N}
\le \sum_{j=1}^{i-1} \vert a_j-b_j\vert + \sum_{j=1}^{i} \vert a_j-b_j\vert.
\end{align*}
Finally, by considering the first and last variates separately, we obtain:
\begin{align*}
\Was_1 (D_{\bf a}, D_{\bf b}) & \le \sum_{i=1}^d \E (\vert X_i-Y_i\vert) \\
& \le 2 (d-1)\vert a_1-b_{1}\vert + 2 (d-2) \vert a_2-b_2\vert + \ldots \\
& \quad + \ldots + 2\times 2 \vert a_{d-2}-b_{d-2}\vert + 2 \vert a_{d-1}-b_{d-1}\vert,
\end{align*}
and this is the claim.
\end{proof}

\subsection{Supplementary figures}
\label{app:figures}

\begin{figure}[ht!]
\centering
\begin{minipage}{0.65\textwidth}
\includegraphics[width=\textwidth]{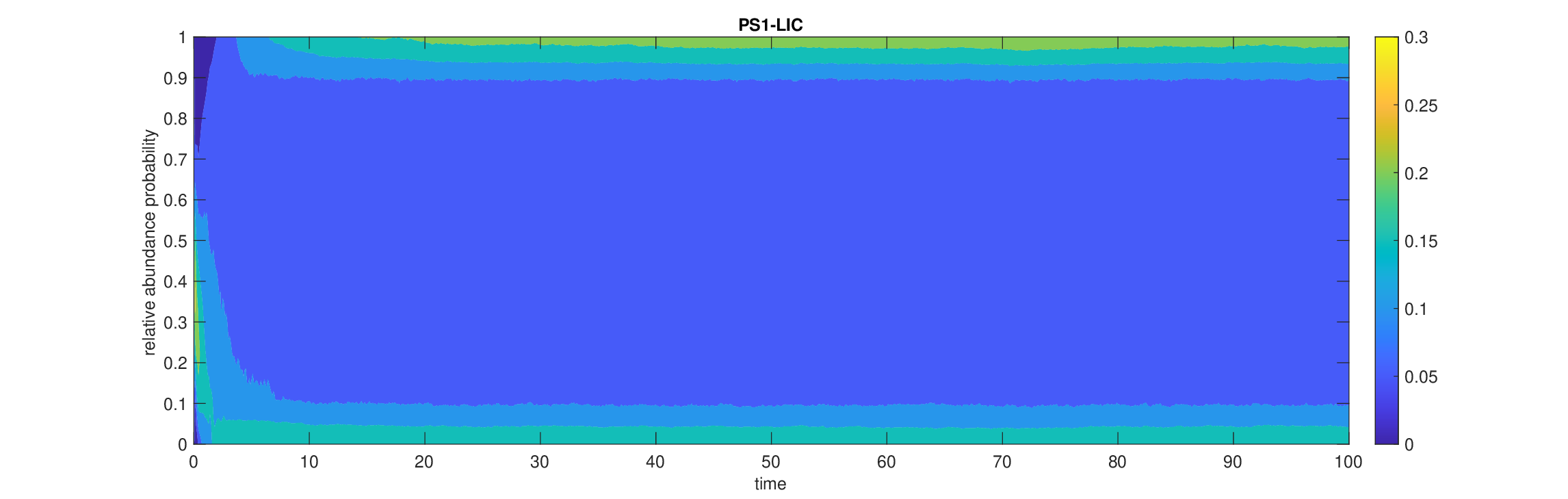}
\end{minipage}
\hfill
\begin{minipage}{0.32\textwidth}
\includegraphics[width=\textwidth]{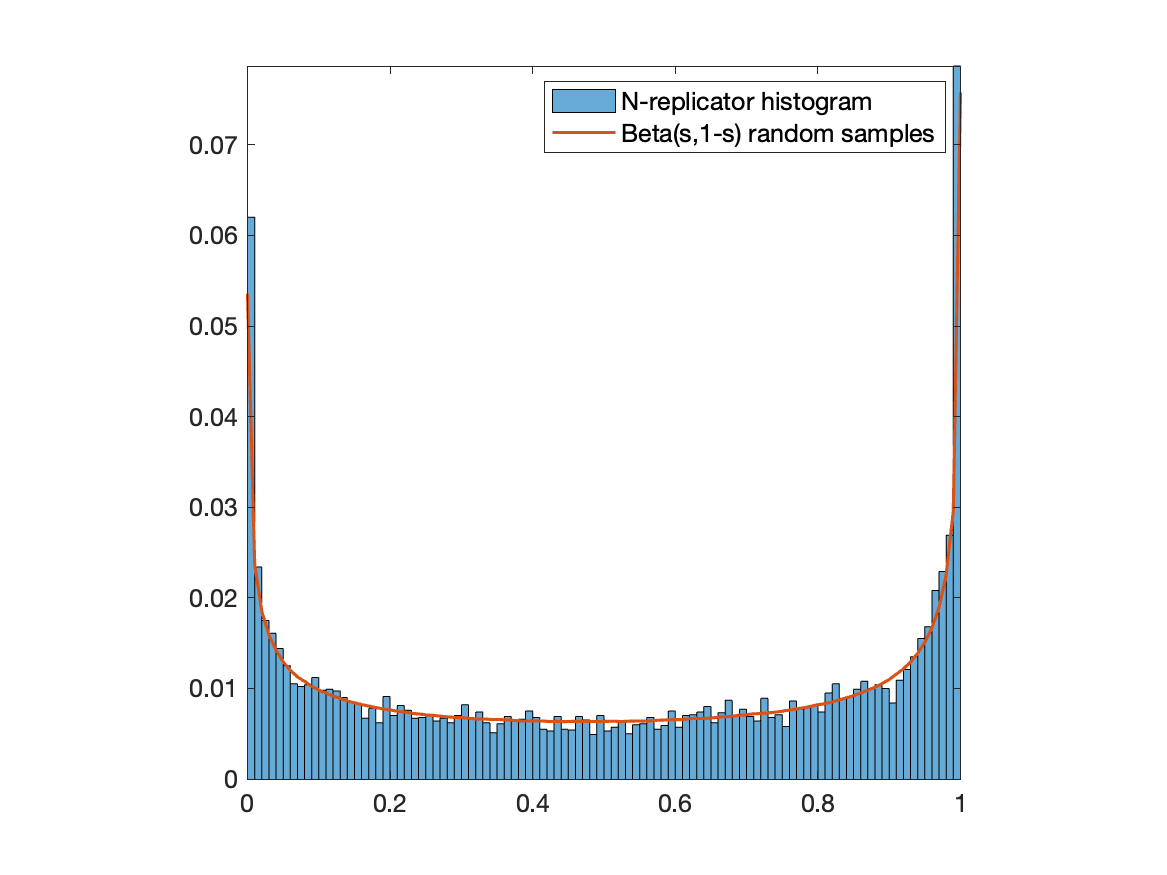}
\end{minipage}

\begin{minipage}{0.65\textwidth}
\includegraphics[width=\textwidth]{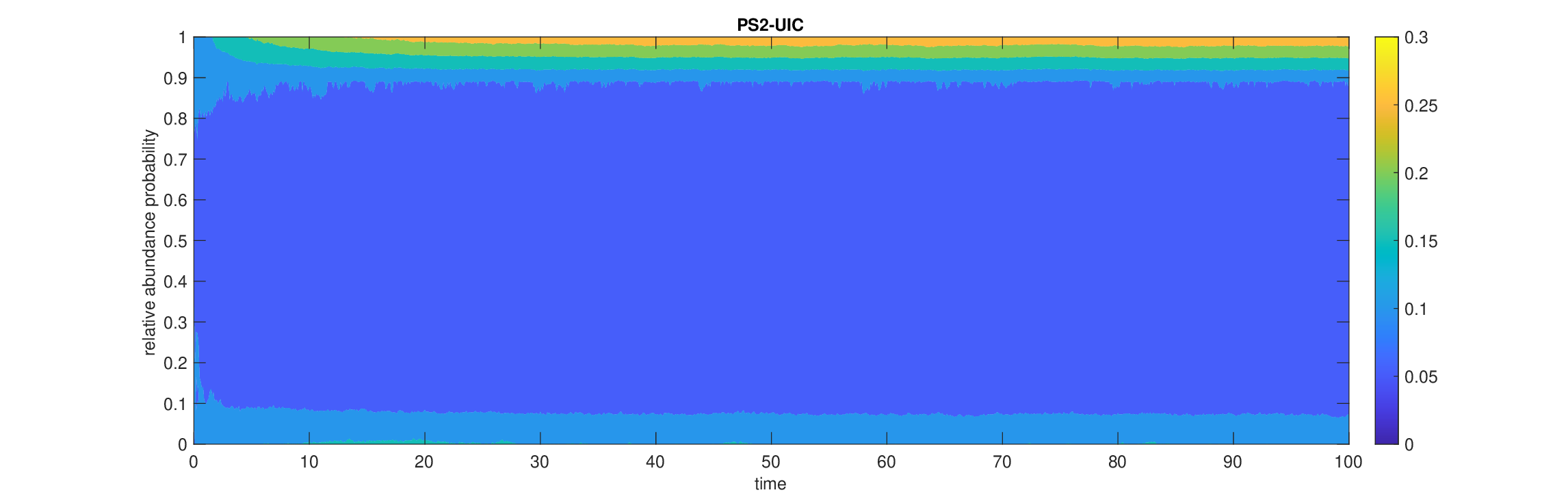}
\end{minipage}
\hfill
\begin{minipage}{0.32\textwidth}
\includegraphics[width=\textwidth]{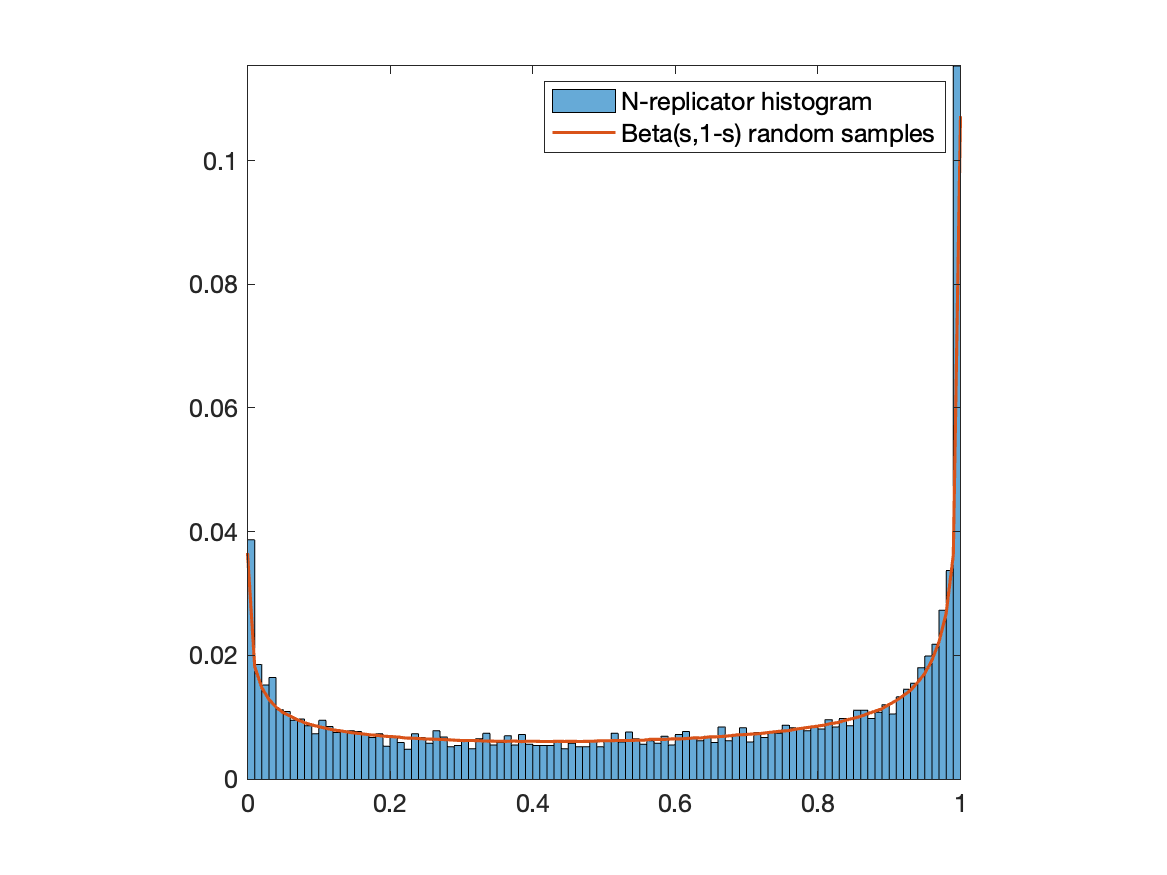}
\end{minipage}

\begin{minipage}{0.65\textwidth}
\includegraphics[width=\textwidth]{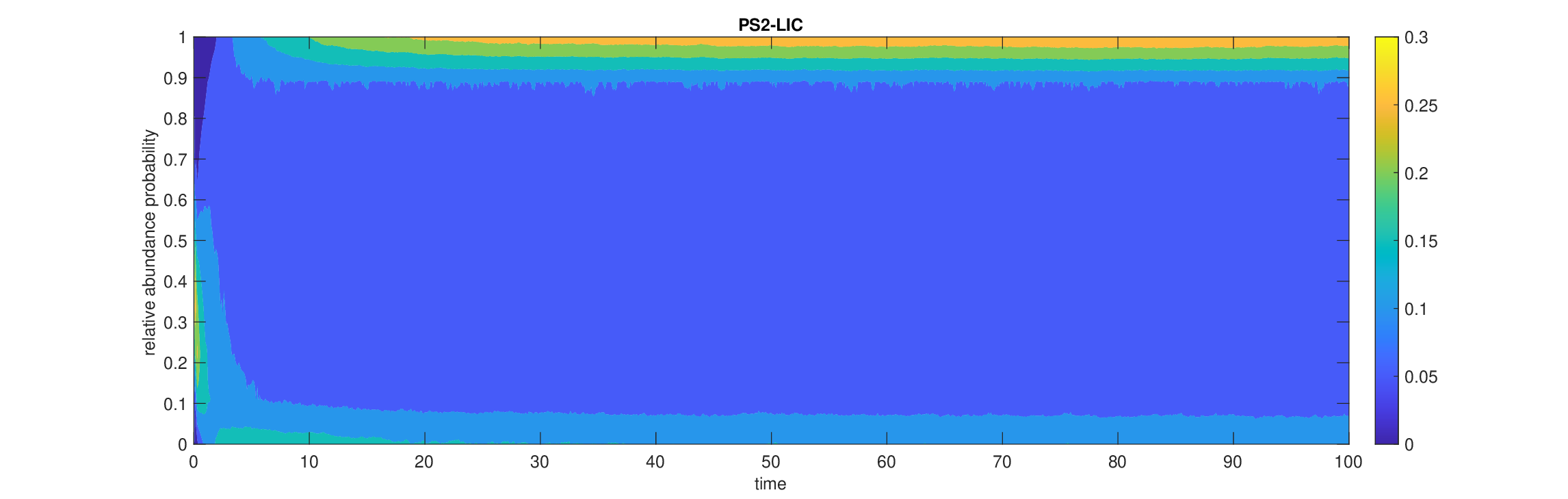}
\end{minipage}
\hfill
\begin{minipage}{0.32\textwidth}
\includegraphics[width=\textwidth]{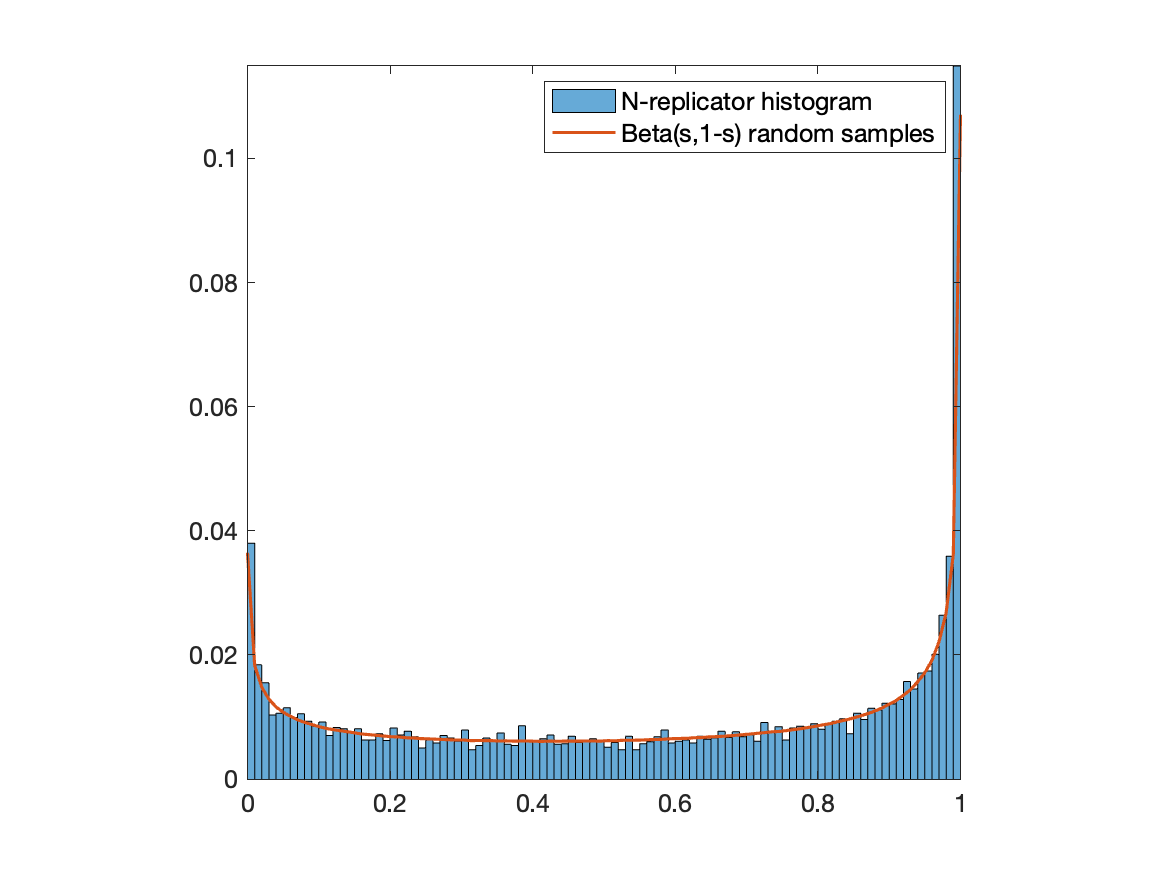}
\end{minipage}

\caption{Simulations for the $(\simp{2})^{10000}$ interacting particle system. For the $N$-replicator system, at each time step, we calculate the empirical histogram from the first coordinate of each replicator under probability normalization. From top to bottom, the plot shows the time evolution of the empirical relative abundance probability for PS1-LIC, PS2-UIC, and PS2-LIC.}
   \label{fig:a2d_simplex}
\end{figure}

In this section, we present some supplementary results for the model studied in Section~\ref{sec:simulations}. Figure~\ref {fig:a2d_simplex}, shows the time evolution of the $N$-replicator for as second set of initial conditions (LIC) and also for the second set of parameters (PS2-UIC and PS2-LIC) along with the comparative between the final histogram for the $N$-replicator dynamics and the Beta theoretical distribution of Theorem~\ref{th:simple_theorem}. These numerical results support the claim that the long-term convergence and the shape of the steady state do not depend on the initial conditions chosen or the used parameter set.

The second supplementary figure shows the result for the goodness-of-fit for different parameter sets at times $t=40$ and $t=50$. This is a graphical representation of the results explained in the main text, along with the 90\%, 95\%, and 99\% quantiles. We can see that in none of the cases we can reject the null hypothesis, supporting the conjecture that as the system evolves in time, the solutions converge to the stationary distribution given by Theorem~\ref{th:simple_theorem}

\begin{figure}[!ht]
\centering
\includegraphics[width=0.45\textwidth]{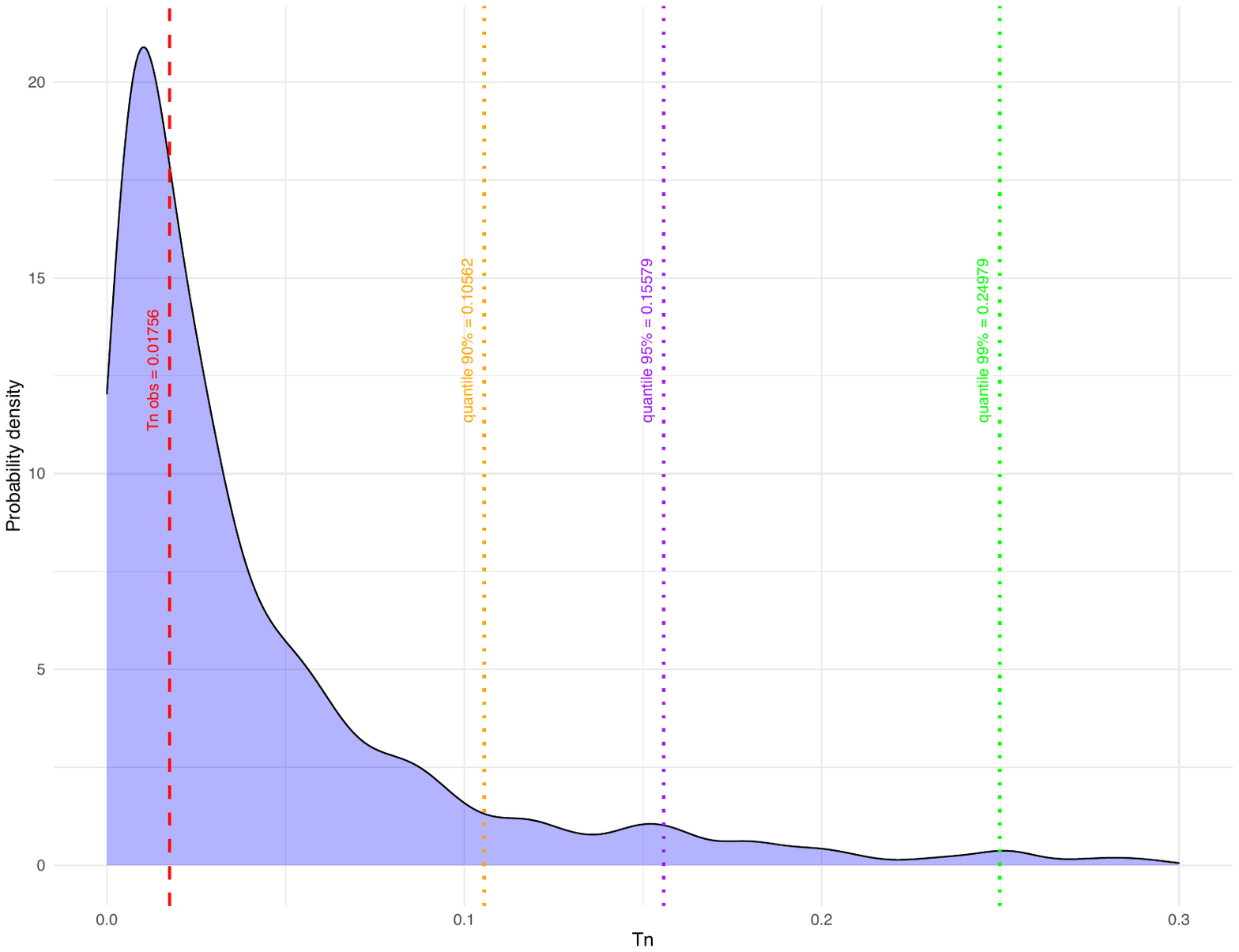}
\includegraphics[width=0.45\textwidth]{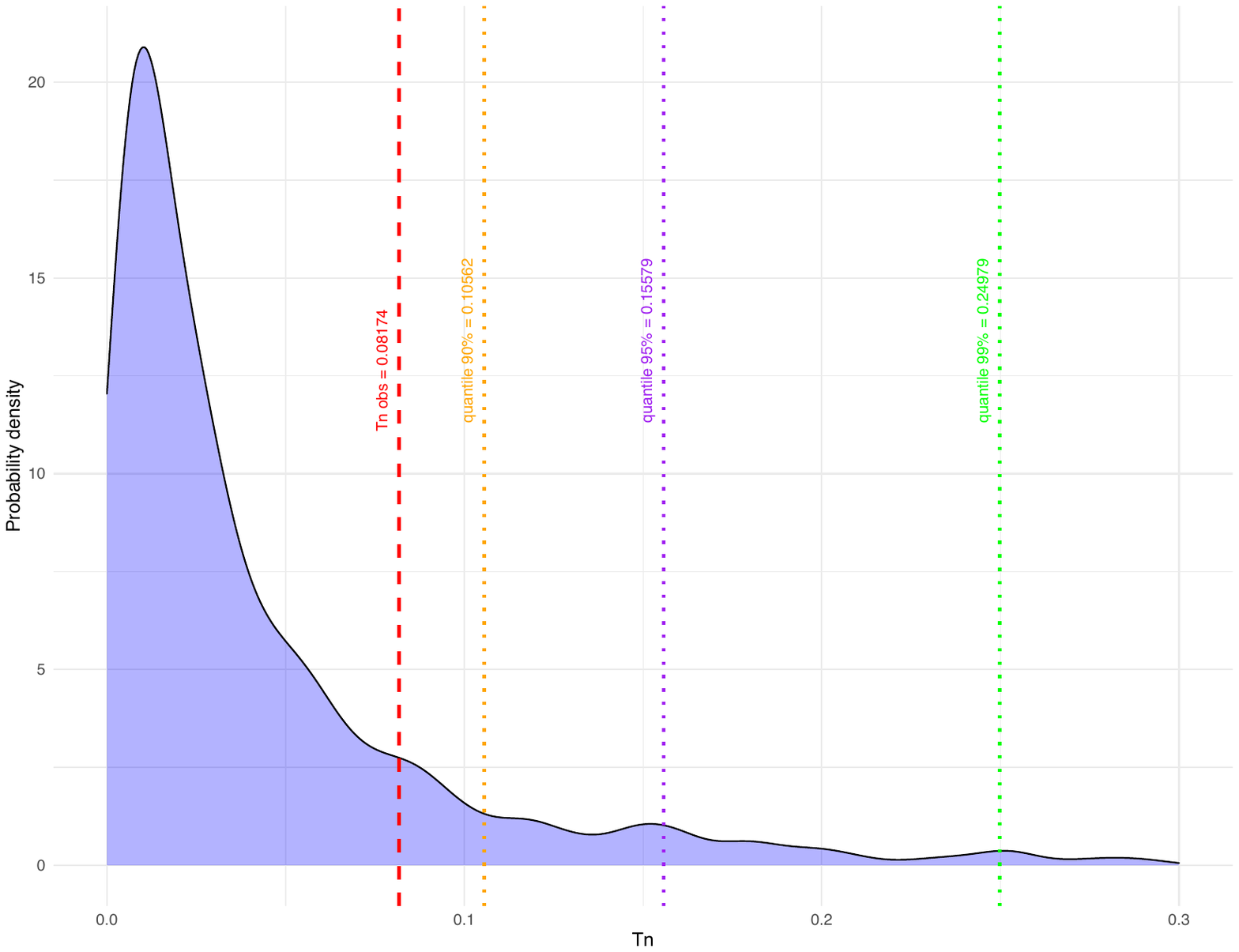}

\includegraphics[width=0.45\textwidth]{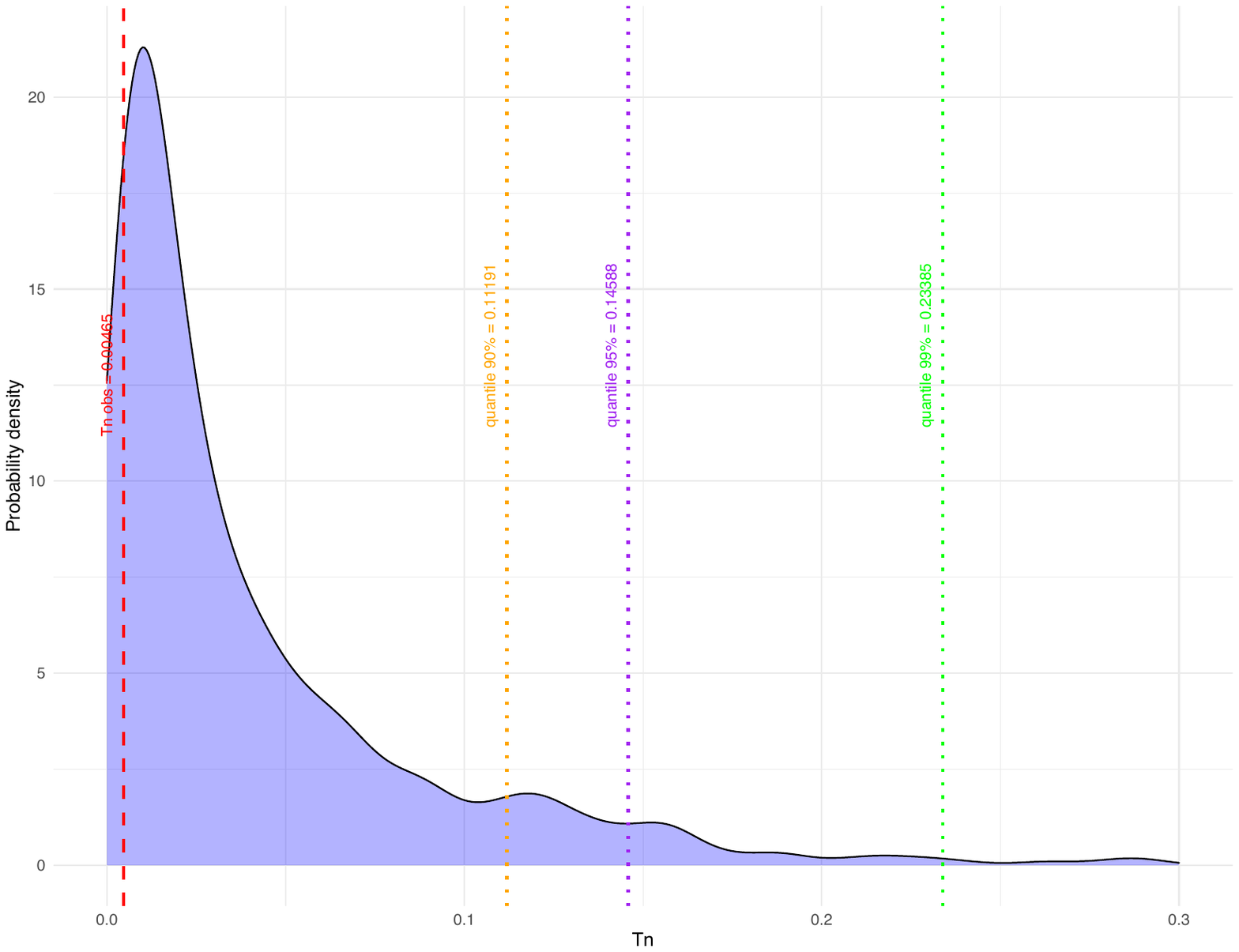}
\includegraphics[width=0.45\textwidth]{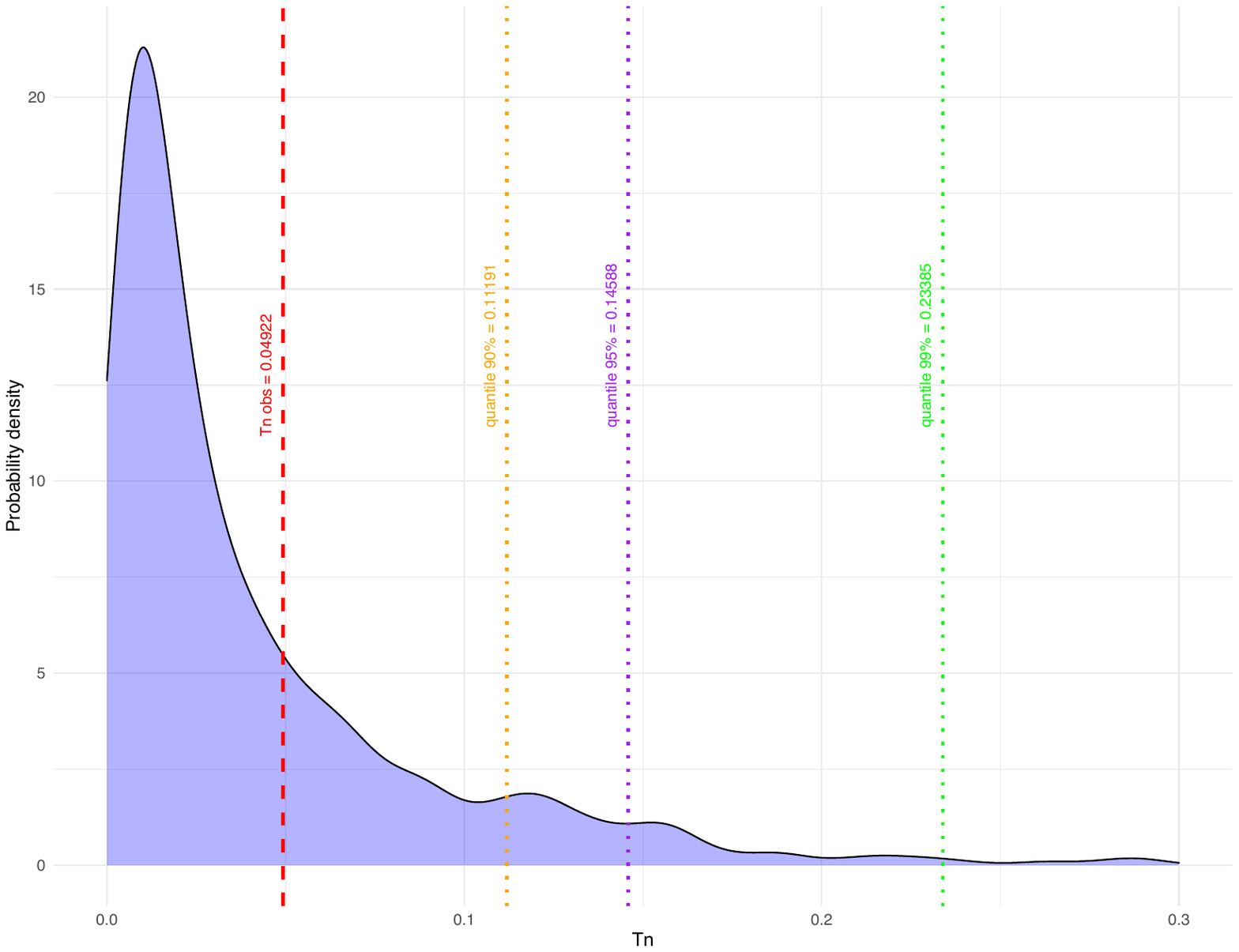}
\caption{Graphical hypothesis testing for contrasting $H_0$: ``the law of the  first coordinate trajectories of~\eqref{eq:stochastic_replicator_2} follows a $Beta(a,b)$ at time $t$''. The statistical test is based on the work of \cite{ebner2021new}. The upper panels show the results for PS1 and the lower panels are the corresponding results for PS2, with their corresponding theoretical $s$ given in the main text. The left column corresponds to time $t=40$ and the right one to time $t=50$. The test is based on the statistic $T_n$, constructed according to a $L^2$ distance given in Equation \eqref{stat-test}, whose observed values are shown by the red dashed lines. We can notice that, in all the cases, they are much lower than the quantiles that delimit the rejecting region at each indicated $(1-\alpha)\times100$ percentile, with $\alpha=0.1$ (yellow dashed lines), $\alpha=0.05$ (purple dashed lines) and  $\alpha=0.01$ (green dashed lines),   regarding the theoretical probability density of the statistic $T_n$ under $H_0$, approximated by simulations. 
}
\label{fig:NewBeta}
\end{figure}

\end{document}